\documentclass[11 pt]{article}
\usepackage{style}
\usepackage{booktabs}
\usepackage{multirow}
\usepackage{graphicx}
\usepackage[caption=false,font=footnotesize]{subfig}

\newenvironment{pf*}[1][]{%
  \par\vspace{1ex}\noindent
  \ifx&#1&%
    \textbf{Proof.}~%
  \else%
    \textbf{#1.}~%
  \fi%
}{\hfill$\qed$\par\vspace{1ex}}
\title{\LARGE\bfseries A Zeroth-Order Extra-Gradient Method for Black-Box Constrained
Optimization}
\author{
Yuke Zhou\textsuperscript{$*,1$}, Ruiyang Jin \textsuperscript{$\dagger,1$}, Siyang Gao\textsuperscript{$\dagger,2$}, Jianxiao Wang\textsuperscript{$*,2$}, and Jie Song\textsuperscript{$*,3$}\\
{\small
\textsuperscript{$*$}\textit{Peking University,} \href{zhouyk@pku.edu.cn}{\textit{\textsuperscript{$1$}zhouyk@pku.edu.cn}},
\href{wang-jx@pku.edu.cn}{\textit{\textsuperscript{$2$}wang-jx@pku.edu.cn}}, 
\href{jie.song@pku.edu.cn}{\textit{\textsuperscript{$3$}jie.song@pku.edu.cn}}}\\
{\small
\textsuperscript{$\dagger$}\textit{City University of Hong Kong,} \href{ruiyajin@cityu.edu.hk}{\textit{\textsuperscript{$1$}ruiyajin@cityu.edu.hk}}, 
\href{siyangao@cityu.edu.hk}{\textit{\textsuperscript{$2$}siyangao@cityu.edu.hk}}}
}
\date{\vspace{-0.4 in}}
\begin{document}
\maketitle

\begin{abstract}
Non-analytical objectives and constraints often arise in control systems, particularly in problems with complex dynamics, which are challenging yet lack efficient solution methods.
In this work, we consider general constrained optimization problems involving black-box objectives and constraints. To solve it, we reformulate it as a min-max problem and propose a zeroth-order extra gradient (ZOEG) algorithm that combines the extra gradient method with a feedback-based stochastic zeroth-order gradient estimator. Then, we apply another coordinate gradient estimator to design the zeroth-order coordinate extra gradient algorithm (ZOCEG) to further improve efficiency. The theoretical analysis shows that ZOEG can achieve the best-known oracle complexity of $\mathcal{O}(d\epsilon^{-2})$ to get an $\epsilon$-optimal solution ($d$ is the dimension of decision space), and ZOCEG can improve it to $\mathcal{O}(d\epsilon^{-1})$. Furthermore, we develop a variant of ZOCEG, which applies block coordinate updates to enhance the efficiency of single-step gradient estimation. Finally, numerical experiments on a load tracking problem validate our theoretical results and the effectiveness of the proposed algorithms.
\end{abstract}

\section{Introduction}
In this paper, we consider a general constrained optimization problem in complex systems where the objective and constraint functions are deterministic but black-box, i.e., we only have access to their zeroth-order feedback (function values) rather than any higher-order information, such as gradients and Hessians. Specifically, we consider the following constrained optimization problem:
\begin{equation}\label{eq: constrained_problem}
    \min_{\mathbf{x}\in\mathcal{X}}\ \phi_0(\mathbf{x})\ \ 
\text {s.t.}\ \phi_j(\mathbf{x})\leq0,\ \forall j=1,2,...,d_y,
\end{equation}
where $\mathcal{X}\in\mathbb{R}^{d_x}$ is a convex and compact set, and $\phi_j:\mathbb{R}^{d_x}\to\mathbb{R},\ \forall j=0,1,...,d_y$ are convex functions. They have no analytical forms, but can be evaluated at any query point $x\in\mathbb{R}^{d_x}$. Such optimization problems with black-box dynamics arise in diverse domains, including optimal control \cite{guo2025safe}, simulation optimization \cite{fu2015handbook}, machine learning \cite{nguyen2023stochastic}, power systems optimization \cite{chen2021modelfreeoptimalvoltagecontrol}, etc. In these problems, the objective and constraints usually lack analytical models due to the high complexities of the dynamics, expensive computation, and privacy issues \cite{jin2023zeroth}.
Therefore, we can only access their zeroth-order feedback via measurements or simulations. 
For instance, we can hardly model the dynamics of a distribution network in power systems when information about the topology and parameters is unavailable. As a result, we can only observe the feedback, such as the voltages and power flow, by performing the given input. 

Generally, traditional gradient-based algorithms, which rely on gradient computation, cannot be directly applied to problems with black-box dynamics. This necessitates the use of derivative-free optimization methods. Among them, direct search approaches are widely adopted \cite{hooke1961direct}, including the simplex methods and directional direct search (DDS) methods. By constructing and improving a simplex according to the evaluated function values at the vertices, the Nelder-Mead simplex algorithm is one of the most popular derivative-free methods \cite{nelder1965simplex}. The DDS methods typically consist of a search step and a poll step. In the search step, a candidate set of points is evaluated to check the improvement in function values. In the poll step, neighboring points are explored when no improvement is detected in the search step. Various algorithms fall under this framework, primarily distinguished by how the poll directions are generated. Specifically, coordinate search uses elementary basis vectors for poll directions \cite{fermi1952numerical}, while generalized pattern search 
(GPS) samples poll directions from a positive spanning set defined on a mesh \cite{torczon1997convergence}, which allows more flexible and adaptive search directions \cite{lewis1999pattern}. Furthermore, the mesh adaptive direct search (MADS) method generalizes GPS to achieve improved performance by adapting the mesh dynamically \cite{audet2006mesh}. Nevertheless, these methods are restricted to low-dimensional problems.
Another line of research focuses on model-based approaches. These studies develop surrogate models to approximate the objective and constraint functions. The commonly used surrogates include polynomial models \cite{conn1991globally}, Gaussian process \cite{ungredda2024bayesian}, Radial Basis Function (RBF) \cite{regis2011stochastic}, etc. 
The algorithm alternates between the following two steps until termination: (i) evaluating the objective values at sampled points and updating the surrogate model; and (ii)  sampling new points based on the surrogate model and past data, a process that often involves solving subproblems \cite{regis2020survey}. However, their performance is highly dependent on the choice of surrogate models. 

Many techniques have been developed in derivative-free optimization to deal with black-box constraints, mainly including penalty, filter, and approximation methods. For penalty methods, a merit function combining the objective function and the constraint violation is minimized \cite{dzahini2023constrained}. Filter methods are inspired by bi-objective optimization, where each update seeks to decrease either the objective function value or the constraint violation \cite{shi2019filter}. For model-based approaches, approximation models are built for both the objective and constraint functions to handle the black-box constraints \cite{regis2017conorbit}. However, most studies mentioned above lack convergence analysis, and the few that establish convergence results only provide asymptotic analysis.

Zeroth-order optimization (ZO), another representative method in derivative-free optimization, has recently gained increasing attention due to its ease of implementation and capability to deal with high-dimensional problems \cite{wang2018stochastic}. 
The core idea behind ZO is to use the finite difference of function values to construct a gradient estimator and embed it into gradient-based algorithms \cite{berahas2022theoretical}. The gradient estimation includes single-point, $2$-point, and $2d$-point approaches, where $d$ is the dimension of the decision variable. Wherein, the $2d$-point coordinate gradient estimator (CooGE) is the most intuitive one that estimates the partial gradient along each dimension asynchronously, but it is inefficient to generate a full gradient for high-dimensional problems \cite {xu2024derivative,liu2025inexact}. The single-point and $2$-point gradient estimators synchronously add perturbations sampled from isotropic distributions, e.g., the Gaussian distribution and uniform distribution on a unit sphere (leading to GauGE and UniGE, respectively), to all dimensions to estimate the full gradients \cite{zhang2024boosting,flaxman2005online,nesterov2017random, tang2023zeroth}. However, $2$-point estimators are usually preferred because the variance of single-point estimators is unbounded, thus leading to poor performance in optimization \cite{duchi2015optimal}. By combining zeroth-order gradient estimators with gradient-based algorithms, e.g., gradient descent and mirror descent, ZO can generally achieve similar convergence rates as their first-order counterparts in unconstrained optimization, differing only by a factor of $d$ \cite{liu2020primer}.

When it comes to constrained optimization, the majority of existing ZO research focuses on problems with black-box objectives but simple and explicit constraint sets. 
These white-box constraints allow efficient Euclidean or Bregman projections or admit linear‐minimization oracles, which enable efficient approaches to maintain feasibility at each iteration. For instance, integrating zeroth-order gradient estimators with projected gradient descent yields the ZO-PSGD algorithm \cite{ghadimi2016mini}, which can be further extended to the mirror descent framework \cite{duchi2015optimal,yu2021distributed}. In addition, projection-free methods have been developed based on the Frank-Wolfe method \cite{chen2020frank,ye2025enhanced}.  
To deal with black-box constraints (as in \eqref{eq: constrained_problem}), primal-dual frameworks are widely employed to solve the derived min-max problems \cite{hu2023gradient}. By combining constraint extrapolation techniques with zeroth-order stochastic gradient estimators, the \textit{oracle complexity} (number of function value queries) of $\mathcal{O}(d^2\epsilon^{-2} )$ is achieved in \cite{nguyen2023stochastic}, where $\epsilon$ is the target error. Besides, the zeroth-order gradient descent ascent (GDA) algorithm was developed to solve the min-max problems \cite{liu2020min}. A variant of zeroth-order GDA, called zeroth-order optimistic GDA, is proposed, and the oracle complexity of $\mathcal{O}(d^4\epsilon^{-2})$ is established \cite{maheshwari2022zeroth}. However, the algorithms tailored to deal with problem \eqref{eq: constrained_problem} are still scarce, and the existing ones suffer from slow convergence due to the large variance of stochastic gradient estimation.

To address the challenges of ZO with black-box constraints, we design a zeroth-order extra-gradient (ZOEG) algorithm that integrates the $2$-point UniGE and a variant of GDA, called extra-gradient (EG) \cite{nemirovski2004prox}, in this work to solve the problem \eqref{eq: constrained_problem}. In addition, we perform finite-sample analysis of ZOEG to achieve improved dependence on the dimension $d$. Furthermore, to address the variances of the $2$-point uniform gradient estimator, we leverage the $2d$-point CooGE with the EG method and design the zeroth-order coordinate extra-gradient (ZOCEG) algorithm to achieve improved oracle complexity. We further extend ZOCEG to zeroth-order block coordinate extra-gradient (ZOBCEG) by randomly selecting subsets of coordinates for updates at each iteration, thereby enhancing the single-step computational efficiency. The main results are summarized in Table \ref{tab: complexity}, and the contributions of this paper are summarized as follows.
\begin{itemize}
    \item We apply the primal-dual framework to solve problem \eqref{eq: constrained_problem} and reformulate it into a min-max problem. Then, we design the ZOEG algorithm by embedding the $2$-point UniGE into the EG method to solve the problem without requiring any higher-order information. The theoretical analysis of ZOEG provides the oracle complexity $O(d\epsilon^{-2})$ to find an $\epsilon$-optimal solution with the best-known dependence on $d$.
    \item To further address the large variance of $2$-point estimators, we improve the proposed ZOEG by applying the coordinate gradient estimators and propose a new algorithm, ZOCEG, which estimates the full gradient information with the $2d$-point CooGE. The theoretical analysis establishes the oracle complexity $O(d\epsilon^{-1})$ of ZOCEG. This result represents, to the best of our knowledge, the state-of-the-art iteration complexity for solving problem \eqref{eq: constrained_problem}. We also extend ZOCEG to ZOBCEG to enhance the efficiency of gradient estimation. 
    \item We validate the proposed algorithms and theoretical results on a practical load tracking task in power systems, and compare our algorithms with existing algorithms. The proposed algorithms show superior convergence and also perform well in nonconvex problems.
\end{itemize}

\begin{table}[htbp]
  \centering
  \caption{Comparison of complexity bounds.}
    \begin{tabular}{|c|c|c|}
    \hline
    \textbf{Algorithm}  & \begin{tabular}[c]{@{}c@{}}\textbf{Gradient}\\ \textbf{Estimator}\end{tabular} & \begin{tabular}[c]{@{}c@{}} \textbf{Oracle}\\ \textbf{Complexity} \end{tabular} \\
   \hline
    OGDA-RR \cite{maheshwari2022zeroth}  & 2-point GauGE & $\mathcal{O}
(d^4\epsilon^{-2})$\\
    \hline
    SZO-ConEX \cite{nguyen2023stochastic} & 2-point GauGE & $\mathcal{O}
    (d^2\epsilon^{-2})$ \\
    \hline
    ZOEG (Ours)  & 2-point UniGE & $\mathcal{O}
    (d\epsilon^{-2})$ \\
   \hline
    ZOCEG (Ours) & $2d$-point CooGE &$\mathcal{O}
    (d\epsilon^{-1})$ \\
    \hline
    ZOBCEG (Ours) & $2\tau$-point CooGE &$\mathcal{O}
    (d\epsilon^{-2})$ \\
    \hline
    \end{tabular}%
  \label{tab: complexity}%
  \vspace{1ex}
  \begin{minipage}{\linewidth}
     \footnotesize{\textit{Note:} $\tau$ denotes the block size (see detailed definition in Section \ref{sec: ZOCEG}).}
  \end{minipage}
\end{table}%

The rest of the paper is organized as follows. Section II introduces some preliminaries of the methodology in our work. In Section III, we present the ZOEG algorithm and establish its finite-sample convergence. The design and convergence results of the ZOCEG algorithm and its extension are presented in Section IV. We verify the performance of the proposed algorithms via numerical experiments in Section V. Finally, Section VI concludes this paper.

\section{Preliminaries}
\textit{Notations.} For an integer $n>0$, we denote $[n]:=1,2,...,n$. For a vector $\mathbf{z}\in\mathbb{R}^d$, we denote the $i$th entry of $\mathbf{z}$ by $z(i),\forall i\in [d]$, and $\mathbf{z}(i:j)$ denotes the subvector of $\mathbf{z}$ consisting of the elements indexed from $i$ to $j$. The notation $[\mathbf{z}]_+$ represents the component-wise maximum of vectors $\mathbf{z}$ and $\mathbf{0}$. For a differentiable function $h(\mathbf{z}):\mathbb{R}^n\to\mathbb{R}$, we denote the partial gradient with respect to 
$z(i)$ by $\nabla_i h(\mathbf{z}),i\in[n]$.
Similarly, given a differentiable function $f(\mathbf{x},\mathbf{y}):\mathbb{R}^{d_x}\times\mathbb{R}^{d_y}\to\mathbb{R}$, we denote the $i$th (and $j$th) coordinate partial gradient in $\mathbf{x}$ (and $\mathbf{y}$) by $\nabla_{\mathbf{x},i}f(\mathbf{x},\mathbf{y})$ (and $\nabla_{\mathbf{y},j}f(\mathbf{x},\mathbf{y})$).
Let $\langle \cdot,\cdot\rangle$ denote the standard inner product and let $\Vert\cdot\Vert$ denote the $\ell_2$ norm in Euclidean space. The unit sphere in $\mathbb{R}^d$ is denoted by $\mathbb{S}_{d-1}:=\{\mathbf{z}\in\mathbb{R}^d:\Vert\mathbf{z}\Vert=1\}$, the $d$-dimensional unit ball in $\mathbb{R}^d$ is denoted by $\mathbb{B}_d:=\{\mathbf{z}\in\mathbb{R}^d:\Vert\mathbf{z}\Vert\leq1\}$, and $\mathcal{U}(\cdot)$ denotes the continuous uniform distribution over a set. For a convex set $\mathcal{Z}$, let $\mathcal{P}_{\mathcal{Z}}$ denote the projection operator to set $\mathcal{Z}$. For any set $\mathcal{S}$, we use 
$|\mathcal{S}|$ to denote its cardinality, i.e., the number of elements in the set.

\subsection{Problem Formulation and Assumptions}
To address the constrained optimization problem \eqref{eq: constrained_problem}, we first apply the primal-dual framework and derive the corresponding Lagrangian function as follows:
$$f(\mathbf{x},\mathbf{y})=\phi_0(\mathbf{x})+\sum_{j=1}^{d_y}y_j\phi_j(\mathbf{x}),$$
which is a convex-concave function, i.e., $f(\mathbf{x},\mathbf{y})$ is convex in $\mathbf{x}$ and concave in $\mathbf{y}$. Then, given the problem \eqref{eq: constrained_problem}, we can instead consider the following general convex-concave min-max problem:
\begin{equation}\label{eq: min-max}
\min_{\mathbf{x}\in\mathcal{X}}\max_{\mathbf{y}\in\mathcal{Y}}f(\mathbf{x},\mathbf{y}),
\end{equation}
where $\mathcal{Y}\subseteq\mathbb{R}^{d_y}$ is a convex constraint set. For notational simplicity, we denote $\mathbf{z}=(\mathbf{x},\mathbf{y})$, and $\mathcal{Z}=\mathcal{X}\times\mathcal{Y}$, $d=d_x+d_y$. Define the operator $F(\mathbf{z}):=[\nabla_{\mathbf{x}}f(\mathbf{x},\mathbf{y});-\nabla_{\mathbf{y}}f(\mathbf{x},\mathbf{y})]$. To evaluate the solution quality of problem \eqref{eq: min-max}, we define the duality gap for any $\mathbf{z}=(\mathbf{x},\mathbf{y})\in \mathcal{Z}$:
$$\Delta(\mathbf{z})=f(\mathbf{x},\mathbf{y}^*)-f(\mathbf{x}^*,\mathbf{y}),$$
where $\mathbf{z}^*=(\mathbf{x}^*,\mathbf{y}^*)$ denotes the min-max saddle point of $f(\mathbf{x},\mathbf{y})$. The existence of $\mathbf{z}^*$ is assumed throughout this research. Accordingly, a feasible point $\mathbf{z}\in\mathcal{Z}$ is said to be an $\epsilon$-saddle point of $f$ if $\Delta(\mathbf{z})\leq \epsilon$. Clearly, $\mathbf{z}$ is the saddle point of $f$ if $\Delta (\mathbf{z})=0$, which provides an optimal solution to problem \eqref{eq: constrained_problem}. Similarly, we say a point $\mathbf{x}\in\mathbb{R}^{d_x}$ is $\epsilon$-optimal for problem \eqref{eq: constrained_problem} if $[\boldsymbol{\phi}(\mathbf{x})]_+\leq \epsilon $ and $\boldsymbol{\phi}_0 (\mathbf{x})-\boldsymbol{\phi}_0^*\leq \epsilon$, where $ \boldsymbol{\phi}=(\phi_1, \dots,\phi_{d_y})^T$ denotes the vector of constraint functions, and $\phi_0^*=\phi_0(\mathbf{x}^*)$ denotes the optimal objective value.

Then, we formally provide some mild assumptions on the problem \eqref{eq: min-max}, which will be used throughout the analysis.
\begin{assumption}\label{ass: convexity}
    The function $f(\mathbf{x},\mathbf{y})$ is continuously differentiable over $\mathcal{X}\times\mathcal{Y}$. Furthermore, $f(\mathbf{x},\mathbf{y})$ is convex-concave, i.e., convex in $x$ and concave in $y$.
\end{assumption}

\begin{assumption}\label{ass: smoothness}
    The function $f$ is $G$-Lipschitz and $L$-smooth, i.e., for any $\mathbf{z}_1,\mathbf{z}_2\in\mathcal{Z}$,
    \begin{equation*}
        \Vert\nabla f(\mathbf{z}_1)-\nabla f(\mathbf{z}_2)\Vert\leq L \Vert\mathbf{z}_1-\mathbf{z}_2\Vert,
    \end{equation*}
    \begin{equation*}
        \Vert f(\mathbf{z}_1)- f(\mathbf{z}_2)\Vert\leq G \Vert\mathbf{z}_1-\mathbf{z}_2\Vert.
    \end{equation*}
\end{assumption}

\begin{assumption}\label{ass: bounded}
    The set $\mathcal{Z}=\mathcal{X}\times\mathcal{Y}$ is compact and convex, thus, we have $\tilde{D}=\max_{\mathbf{z},\mathbf{z}'\in\mathcal{Z}}\Vert\mathbf{z}-\mathbf{z}'\Vert\leq +\infty$ .
\end{assumption}

Assumption \ref{ass: convexity} can be equivalently imposed on the primal problem, i.e., $\boldsymbol{\phi}_j(\mathbf{x}),\forall j\in\{0,1,\cdots,d_y\}$ is continuously differentiable and convex. Assumption \ref{ass: smoothness} imposes the Lipschitz continuity on $f(\mathbf{x},\mathbf{y})$ and its gradients, which is common in smooth optimization \cite{liu2020min, jin2023zeroth}. As for Assumption \ref{ass: bounded}, the boundedness of $\mathcal{X}$ is natural because $x$ is generally bounded in practical problems.
Although $\mathbf{y}\geq 0$ is the Lagrange multiplier and unbounded, the boundedness of $\mathbf{y}^*$ is justified in \cite{nedic2009subgradient} under the Slater condition. Therefore, we can construct a bounded set $\mathcal{Y}$ containing $\mathbf{y}^*$ to replace $\{\mathbf{y}\in \mathbb{R}^{d_y}|\mathbf{y}\geq 0\}$.



\subsection{Zeroth-Order Gradient Estimators}
Note that we can only obtain the function values of $\phi_0(\mathbf{x})$ and $\phi_j(\mathbf{x})$ for all $j \in [d_y]$, but lack access to higher-order information such as gradients and Hessians. This limitation renders traditional first-order min-max optimization methods inapplicable. To address this, we turn to ZO methods by estimating gradients using function values. 
An intuitive approach for estimating the full gradient involves computing partial derivatives along coordinate directions using finite differences.
Given a differentiable function $h(\mathbf{z}):\mathbb{R}^{d}\to\mathbb{R}$, the partial gradient along the $i$th coordinate direction is estimated using the finite difference of two function values:
\begin{equation}  \label{eq: coor_grad_est}
G_{h,i}(\mathbf{z};r)=\frac{h(\mathbf{z}+r\cdot\mathbf{e}_i)-h(\mathbf{z})}{r},
\end{equation}
where $\mathbf{e}_i\in\mathbb{R}^{d}$ is the unit vector with only the $i$th entry being 1. The parameter $r>0$, known as the smoothing radius, controls the magnitude of the perturbation to $\mathbf{z}$. 
By evaluating along each of the $d$ coordinate directions, we obtain the $2d$-point CooGE:
\begin{equation*}
   G_h(\mathbf{z};r)=(G_{h,1}(\mathbf{z};r),G_{h,2}(\mathbf{z};r),\cdots,G_{h,d}(\mathbf{z};r))^T. 
\end{equation*}
However, as the dimension $d$ increases, the estimator's requirement of $2d$ function evaluations may lead to high computational overhead for each single step.

On the other hand, the widely applied $2$-point UniGE estimates a full gradient using only two points \cite{malik2020derivative}:
\begin{equation} \label{eq: 2p_grad_est}G_h(\mathbf{z};r;\mathbf{v})=\frac{h(\mathbf{z}+r\cdot\mathbf{v})-h(\mathbf{z})}{r/d}\cdot\mathbf{v},
\end{equation}
where $\mathbf{v}$ is a perturbation vector uniformly sampled on a unit sphere $\mathcal{U}(\mathbb{S}_{d-1})$. 



Finally, we provide two lemmas to characterize the bias and variance of the coordinate-wise partial gradient estimator \eqref{eq: coor_grad_est} and $2$-point UniGE \eqref{eq: 2p_grad_est} under proper assumptions.
\begin{lemma}\label{lem: coor_grad_estimator}
    Assume a continuously differentiable function $h$ is $L$-smooth, then, the bias of \eqref{eq: coor_grad_est} satisfies:
    \begin{align}\label{eq: coor_grad_est_bias}
        \left|G_{h,i}(\mathbf{z};r)-\nabla_i h(\mathbf{z})\right|\leq \frac{1}{2}Lr.
    \end{align}
\end{lemma}

\begin{lemma}\label{lem: 2P_grad_estimator}
Assume a continuously differentiable function $h$ is $G$-Lipschitz and $L$-smooth. Define the smoothed function as $h^r(\mathbf{z})=\mathbb{E}_{\Bar{\mathbf{v}}\sim\mathcal{U}(\mathbb{B}_d)}\left[h(\mathbf{z}+r\Bar{\mathbf{v}})\right]$.
    Then, the estimator \eqref{eq: 2p_grad_est}: satisfies:
    \begin{align}
        \nabla h^r(\mathbf{z})=\mathbb{E}_{\mathbf{v}}\left[G_h(\mathbf{z};r;\mathbf{v})\right],\label{eq: exp_2p_grad}\\
        \Vert \nabla h^r(\mathbf{z})-\nabla h(\mathbf{z})\Vert \leq Lr,\label{eq: diff smooth}
    \end{align}
    Furthermore, the variance of \eqref{eq: 2p_grad_est} is bounded by:
    \begin{equation}\label{eq: var_2p_grad}
        \mathbb{E}_{\mathbf{v}}\left[\Vert G_h(\mathbf{z};r;\mathbf{v})-\nabla h^r(\mathbf{z})\Vert^2\right]\leq \frac{9d^2G^2}{d+2}+\frac{3d^2 L^2 r^2}{4}.
    \end{equation}   
\end{lemma}
The proof of Lemma \ref{lem: coor_grad_estimator} is deferred to Appendix \ref{apx: pf_coor_grad_est}. Lemma~\ref{lem: 2P_grad_estimator} follows from the continuity and smoothness of $h$, with \eqref{eq: exp_2p_grad}-\eqref{eq: diff smooth} provided in \cite{malik2020derivative} and \eqref{eq: var_2p_grad} provided in \cite{berahas2022theoretical}. Lemmas \ref{lem: coor_grad_estimator} and \ref{lem: 2P_grad_estimator} suggest that the biases of both estimators can be small by adopting a small $r$, but the coordinate-wise partial gradient estimator enjoys a better variance bound of $\frac{L^2r^2}{4}$ when $r$ is small.

\section{Zeroth-Order Extra-Gradient}\label{sec: ZOEG}
In this section, we leverage the $2$-point UniGE to design the zeroth-order extra gradient (ZOEG) algorithm to solve problem \eqref{eq: min-max}, and provide its performance guarantees.
\subsection{Algorithm Design}
To solve the min-max problem \eqref{eq: min-max}, we first propose the ZOEG algorithm based on \eqref{eq: 2p_grad_est} and the EG method, as detailed in Algorithm~\ref{alg: 2PZOEG}. The iterations in our algorithm follow the main steps of extra-gradient as follows:
\begin{equation*}
\mathbf{z}^+_k=\mathcal{P}_{\mathcal{Z}}\left [\mathbf{z}_k-\eta F(\mathbf{z}_k)\right],
\end{equation*}
\begin{equation*}
\mathbf{z}_{k+1}=\mathcal{P}_{\mathcal{Z}}\left [\mathbf{z}_k-\eta F(\mathbf{z}_k^+)\right],
\end{equation*}
where $\eta>0$ is the step size and $\mathbf{z}^+_k$ is the mid-point to approximate the proximal point (PP) update. Recall that the PP update: $\mathbf{z}_{k+1}=\mathcal{P}_{\mathcal{Z}}\left [\mathbf{z}_k-\eta F(\mathbf{z}_{k+1})\right]$ relies on gradient evaluation at $\mathbf{z}_{k+1}$\cite{nemirovski2004prox}.
The implicit nature of the PP update makes it computationally intractable in practice. Therefore, the EG algorithm uses the gradient information at the intermediate point $\mathbf{z}_k^+$ to approximate the gradient at $\mathbf{z}_{k+1}$, which effectively generates a proximal direction \cite{mokhtari2020unified}.

Given the absence of derivative information, for the $k$th iterate, we sample $\mathbf{v}$ and $\mathbf{w}$ uniformly from the unit sphere $\mathbb{S}_{d-1}$ to construct $\mathbf{g}_k$ and $\mathbf{g}_k^+$ as the estimations of $F(\mathbf{z}_k)$ and $F(\mathbf{z}_k^+)$ respectively:
\begin{align*}
\mathbf{g}_k=\frac{f(\mathbf{z}_k+r_k\mathbf{v}_k)-f(\mathbf{z}_k)}{r_k/d}\tilde{\mathbf{v}}_k,\\
\mathbf{g}_k^+=\frac{f(\mathbf{z}_k^++r_k\mathbf{w}_k)-f(\mathbf{z}_k^+)}{r_k/d}\tilde{\mathbf{w}}_k,
\end{align*}
where $\tilde{\mathbf{v}}_k=[\mathbf{v}_k(1:d_x);-\mathbf{v}_k(d_x+1:d)]$, and $\tilde{\mathbf{w}}_k=[\mathbf{w}_k(1:d_x);-\mathbf{w}_k(d_x+1:d)]$. According to Lemma \ref{lem: 2P_grad_estimator}, it can be shown that they are good estimators of $F(\mathbf{z}_k)$ and $F(\mathbf{z}_k^+)$ with bounded bias under Assumption \ref{ass: smoothness}, so we replace the first-order derivatives in extra-gradient method with $\mathbf{g}_k$ and $\mathbf{g}_k^+$ respectively \cite{malik2020derivative}.

Note that we sample $\mathbf{v}_k$ (and $\mathbf{w}_k$) from $\mathcal{U}(\mathbb{S}_{d-1})$ to estimate $F(\mathbf{z}_k)$ (and $F(\mathbf{z}_k^+)$) jointly. we can also estimate the partial gradients in $\mathbf{x}$ and $\mathbf{y}$ separately by adding two perturbations from $\mathcal{U}(\mathbb{S}_{d_x-1})$ and $\mathcal{U}(\mathbb{S}_{d_y-1})$ to $\mathbf{x}$ and $\mathbf{y}$ respectively. Moreover, the ZOEG algorithm outputs the averaged value over all iterations, which is only for the convenience of theoretical analysis. In practice, we can take the last iterate as the final decision.


\begin{algorithm}[htpb]
	\caption{Zeroth-Order Extra-Gradient (ZOEG)}
\begin{algorithmic}\label{alg: 2PZOEG}
	\STATE \textbf{Inputs:} Number of iterations $K$, step sizes $\eta$, smoothing radii $\{r_k\}_{k=0}^{K-1}$, and initial point $\mathbf{z}_0$.
	\FOR{$k\xleftarrow{}0$ to $K-1$} 
		\STATE Sample $\mathbf{v}_k$ and $\mathbf{w}_k$ from $\mathcal{U}(\mathbb{S}_{d-1})$ independently.

        \STATE Update $\mathbf{z}_k$ by:
        \begin{equation}  \label{eq: ZOEG_plus}
        \mathbf{z}^+_k=\mathcal{P}_{\mathcal{Z}}\left [\mathbf{z}_k-\eta\mathbf{g}_k\right],
        \end{equation}
        \begin{equation}\label{eq: ZOEG_next}
        \mathbf{z}_{k+1}=\mathcal{P}_{\mathcal{Z}}\left [\mathbf{z}_k-\eta\mathbf{g}_k^+\right].
        \end{equation}
	\ENDFOR
    \STATE \textbf{Outputs:} $\hat{\mathbf{z}}_K=\frac{1}{K}\sum_{k=0}^{K-1} \mathbf{z}^+_k .$
\end{algorithmic}
\end{algorithm}
\subsection{Finite-Sample Analysis}
The following theorem provides the convergence results for the proposed ZOEG algorithm. 
\begin{theorem}\label{thm: ZOEG_duality gap}
 Suppose Assumptions \ref{ass: convexity}-\ref{ass: bounded} hold. The sequences $\{\mathbf{z}_k\}$, and $\{\mathbf{z}_k^+\}$ are the iterates generated by ZOEG Algorithm. Set the smoothing radii to satisfy 
$\sum_{k=1}^{\infty}r_k=M_1<\infty$ and $\sum_{k=1}^{\infty}r_k^2=M_2<\infty$. Then for all $K>1$, the duality gap of the averaged iterates $\hat{\mathbf{z}}_K=\frac{1}{K}\sum_{k=0}^{K-1}\mathbf{z}_k$ satisfies
    \begin{equation*}
\mathbb{E}\left[\Delta(\hat{\mathbf{z}}_K)\right]\leq\frac{\tilde{D}^2}{2\eta K}+\frac{C_1 \eta + C_2}{K}+108\eta d G^2,
    \end{equation*}
    where  $C_1=(9d^2+12)L^2M_2$, $C_2=L\tilde{D}M_1$.
\end{theorem}
\begin{pf*}[Proof of Theorem \ref{thm: ZOEG_duality gap}]
We start by defining the following  filtration $\mathcal{F}_k:=\sigma(\mathbf{z}_0,\mathbf{v}_0,\mathbf{w}_0,\mathbf{z}_1,\mathbf{v}_1,\mathbf{w}_1,...,\mathbf{z}_k).$ Then, the upper bound on the conditional expectation of zeroth-order EG iterates $\left\langle\eta F(\mathbf{z}_k^+),\mathbf{z}_k^+-\mathbf{z}\right\rangle$ can be established in the following lemma, the proof of which is given in Appendix~\ref{pf: bound_EG_term}.

\begin{lemma}\label{lem: bound_EG_term}
    Under Assumption~\ref{ass: convexity} -~\ref{ass: bounded}, given the zeroth-order extra-gradient update rules \eqref{eq: ZOEG_plus}, \eqref{eq: ZOEG_next}, we have for any $\mathbf{z}\in\mathcal{Z}$:
\begin{equation}    
\label{eq: lem cond ZOEG}
\ \mathbb{E}\left[\left\langle\eta F(\mathbf{z}_k^+),\mathbf{z}_k^+-\mathbf{z}\right\rangle\mid \mathcal{F}_k\right]
    \leq \mathbb{E}\left[H_{k}(\mathbf{z})-H_{k+1}(\mathbf{z})+\Delta_k\mid \mathcal{F}_k\right],
\end{equation}
where $H_{k}(\mathbf{z})$ and $\Delta_k$ are defined as $H_{k}(\mathbf{z})=\frac{1}{2}\Vert \mathbf{z}_k-\mathbf{z}\Vert ^2,
\Delta_k=\!\left(\frac{108d^2}{d+2}+12\right)\eta^2G^2\!+\eta^2L^2 (9d^2\!+\!12) r_k^2\!+\eta L \tilde{D}r_k.$
\end{lemma}

Given Lemma~\ref{lem: bound_EG_term}, taking the total expectation of \eqref{eq: lem cond ZOEG} leads to:
\begin{equation*}
     \mathbb{E}\left[\left\langle\eta F(\mathbf{z}_k^+),\mathbf{z}_k^+-\mathbf{z}\right\rangle\right]
     \leq \mathbb{E}\left[H_{k}(\mathbf{z})-H_{k+1}(\mathbf{z})\right]+\Delta_k.
\end{equation*}
Then, by calculating the telescoping sum of \begin{math}
\mathbb{E}[\langle F(\mathbf{z}_k^+),\allowbreak\mathbf{z}_k^+-\mathbf{z}\rangle]\end{math} for $k=0,1,\cdots,K-1$, we can bound $\sum_{k=0}^{K-1}\mathbb{E}[\langle F(\mathbf{z}_k^+),\allowbreak\mathbf{z}_k^+-\mathbf{z}\rangle]$ as follows:
\begin{align*}
    &\ \sum_{k=0}^{K-1}\mathbb{E}\left[\left\langle F(\mathbf{z}_k^+),\mathbf{z}_k^+-\mathbf{z}\right\rangle\right]\\
    \leq&\ \sum_{k=0}^{K-1}\mathbb{E}\left[\frac{H_{k}(\mathbf{z})}{\eta}-\frac{H_{k+1}(\mathbf{z})}{\eta}\right]+\sum_{k=0}^{K-1}\frac{\Delta_k}{\eta}\\
    =&\ \mathbb{E}\left[\frac{H_{0}(\mathbf{z})}{\eta}-\frac{H_{K}(\mathbf{z})}{\eta}\right]+\sum_{k=0}^{K-1}\frac{\Delta_k}{\eta}\\
    \leq&\ \frac{1}{\eta}\mathbb{E}\left[\frac{1}{2}\Vert \mathbf{z}_{0}-\mathbf{z}\Vert ^2\right]+\sum_{k=0}^{K-1}\frac{\Delta_k}{\eta}\\
    \leq&\ \frac{\tilde{D}^2}{2\eta}+\sum_{k=0}^{K-1}\frac{\Delta_k}{\eta},
\end{align*}
where the second inequality follows from the definition of $H_k(\mathbf{z})$ and the fact $\frac{1}{2}\Vert\mathbf{z}_K-\mathbf{z}\Vert^2\geq 0$, and the last inequality comes from the boundedness of $\mathcal{Z}$.
Consequently, we have: 
\begin{equation}\label{eq: bound_ZOEG_telescope}
\sum_{k=0}^{K-1}\mathbb{E}\left[\left\langle F(\mathbf{z}_k^+),\mathbf{z}_k^+-\mathbf{z}\right\rangle\right]
    \leq \frac{\tilde{D}^2}{2\eta}+C_1 \eta + C_2+108K\eta d G^2.
\end{equation}

Then, the relation between the telescoping sum $\sum_{k=0}^{K-1}\mathbb{E}\left[\left\langle F(\mathbf{z}_k^+),\mathbf{z}_k^+-\mathbf{z}\right\rangle\right]$ and the duality gap is given in the following lemma.
\begin{lemma}\label{lem: duality gap}
Suppose Assumptions \ref{ass: convexity}-\ref{ass: bounded} hold. When there exists some bound $\Gamma_K\geq 0$ satisfying $\sum_{k=0}^{K-1}\mathbb{E}[\langle F(\mathbf{z}_k^+),\allowbreak\mathbf{z}_k^+-\mathbf{z}\rangle]
\leq \Gamma_K$ for all $K>1$, $\mathbf{z}\in\mathcal{Z}$, the duality gap of the averaged iterate $\hat{\mathbf{z}}_K$ is bounded as
\begin{align*}
\mathbb{E}\left[\Delta(\hat{\mathbf{z}}_K)\right]\leq \frac{\Gamma_K}{K}.  
\end{align*}
\end{lemma}

\begin{pf*}[Proof of Lemma \ref{lem: duality gap}]
Following the definition of $(\hat{\mathbf{x}}_K,\hat{\mathbf{y}}_K)$, the duality gap $\mathbb{E}\left[f\left(\hat{\mathbf{x}}_K,\mathbf{y}^*\right)-f\left(\mathbf{x}^*,\hat{\mathbf{y}}_K\right)\right]$ is upper bounded by:
\begin{align*}
&\mathbb{E}\left[f\left(\hat{\mathbf{x}}_K,\mathbf{y}^*\right)-f\left(\mathbf{x}^*,\hat{\mathbf{y}}_K\right)\right]\\
=&\mathbb{E}\left[f\left(\frac{1}{K}\sum_{k=0}^{K-1}\mathbf{x}_k^+,\mathbf{y}^*\right)-f\left(\mathbf{x}^*,\frac{1}{K}\sum_{k=0}^{K-1}\mathbf{y}_k^+\right)\right]\\
\leq& \frac{1}{K}\sum_{k=0}^{K-1}\mathbb{E}\left[f(\mathbf{x}_k^+,\mathbf{y}^*)-f(\mathbf{x}^*,\mathbf{y}_k^+)\right]\\
=& \frac{1}{K}\!\sum_{k=0}^{K-1}\!\mathbb{E}\left[f(\mathbf{x}_k^+,\mathbf{y}^*)\!-\! f(\mathbf{z}_k^+)\!+\!f(\mathbf{z}_k^+)\!-\!f(\mathbf{x}^*,\mathbf{y}_k^+)\right]\\
\leq&\frac{1}{K}\sum_{k=0}^{K-1}\mathbb{E}\left[\left\langle F(\mathbf{z}_k^+),\mathbf{z}_k^+-\mathbf{z}^*\right\rangle\right]\\
\leq& \frac{\Gamma_K}{K},
\end{align*}
where the first inequality follows from Jensen's inequality and the second one follows from that $f$ is convex in $\mathbf{x}$ and concave in $\mathbf{y}$.
\end{pf*}

Finally, combining the results in \eqref{eq: bound_ZOEG_telescope} and Lemma~\ref{lem: duality gap} finishes the proof.
\end{pf*}

Theorem~\ref {thm: ZOEG_duality gap} shows that, under constant step sizes, the expected duality gap decays at the rate of $\mathcal{O}(1/K)$ up to a fixed error, which can approach zero by employing diminishing or sufficiently small step sizes. 

The theoretical guarantee in Theorem \ref{thm: ZOEG_duality gap} only characterizes the finite-iteration bound of the duality gap of problem \eqref{eq: min-max}, while the behavior of the averaged iterate $\hat{\mathbf{x}}_k=\frac{1
}{K}\sum_{k=1}^K\mathbf{x}_k^+$ w.r.t the original constrained problem \eqref{eq: constrained_problem} remains unclear. To bridge this gap, we present the following proposition, which gives explicit bounds on the optimality gap and constraint violation of $\hat{\mathbf{z}}_K$ for problem \eqref{eq: constrained_problem}.

\begin{proposition}\label{prop: ZOEG vio obj}
Suppose the assumptions and conditions in Theorem \ref{thm: ZOEG_duality gap} hold. Then for any $K>1$, we have
\begin{align*}
\Vert\left[\mathbb{E}\left[ \boldsymbol{\phi}(\hat{\mathbf{x}}_K)\right]\right]_+\Vert &\leq \frac{\tilde{D}^2}{2\eta RK}+\frac{C_1 \eta + C_2}{RK}+\frac{108\eta d G^2}{R},\\
\mathbb{E}\left[ \phi_0(\hat{\mathbf{x}}_K)\right]-\phi_0^* & \leq \frac{\tilde{D}^2}{2\eta K}+\frac{C_1 \eta + C_2}{K}+108\eta d G^2.
\end{align*}
\end{proposition}

\begin{pf*} [Proof of Proposition \ref{prop: ZOEG vio obj}]
    Recall that Lemma~\ref{lem: duality gap} establishes the upper bound for the duality gap in solving the min-max formulation. Following a similar idea, we provide the following lemma to link the telescoping sum $\sum_{k=0}^{K-1}\mathbb{E}\left[\left\langle F(\mathbf{z}_k^+),\mathbf{z}_k^+-\mathbf{z}\right\rangle\right]$ to both the performance in the optimality and constaint violation of the derived solution $\hat{\mathbf{x}}_K$. 
\begin{lemma}\label{lem: optimality gap}
Suppose Assumptions \ref{ass: convexity}-\ref{ass: bounded} hold. When there exists some bound $\Gamma_K\geq 0$ satisfying $\sum_{k=0}^{K-1}\mathbb{E}[\langle F(\mathbf{z}_k^+),\allowbreak\mathbf{z}_k^+-\mathbf{z}\rangle]\leq \Gamma_K$ for any $K>1$, and $\mathbf{z}\in\mathcal{Z}$, the expected optimality gap and constraint violation are bounded by
\begin{align*}
\left\Vert\left[\mathbb{E}\left[ \boldsymbol{\phi}(\hat{\mathbf{x}}_K)\right]\right]_+\right\Vert &\leq \frac{\Gamma_K}{RK},\\
\mathbb{E}\left[ \phi_0(\hat{\mathbf{x}}_K)\right]-\phi_0^* &\leq \frac{\Gamma_K}{K}, 
\end{align*}
where $R=\max \{r>0|\mathbf{y}^*+r\mathbf{u}\in\mathcal{Y}, \Vert\mathbf{u}\Vert=1,\forall \mathbf{u}\in\mathbb{R}^{d_y}\}$.
\end{lemma}
\begin{pf*}[Proof of Lemma \ref{lem: optimality gap}]
We start from the upper bound on $\left\langle \mathbf{y}-\mathbf{y}^*,\nabla f_{\mathbf{y}}(\mathbf{z}_k^+)\right\rangle, \forall \mathbf{y}\in\mathcal{Y}$, as follows:
\begin{align}\label{eq: cons_violation}
&\ \left\langle \mathbf{y}-\mathbf{y}^*,\nabla f_{\mathbf{y}}(\mathbf{z}_k^+)\right\rangle \nonumber\\
= &\ \left\langle \mathbf{y}-\mathbf{y}_k^+,\nabla f_{\mathbf{y}}(\mathbf{z}_k^+)\right\rangle + \left\langle \mathbf{y}_k^+-\mathbf{y}^*,\nabla f_{\mathbf{y}}(\mathbf{z}_k^+)\right\rangle\nonumber\\
\leq &\ \left\langle \mathbf{y}_k^+ - \mathbf{y},-\nabla f_{\mathbf{y}}(\mathbf{z}_k^+)\right\rangle +f(\mathbf{z}_k^+)-f(\mathbf{x}_k^+,\mathbf{y}^*)\nonumber\\
\leq &\ \left\langle \mathbf{y}_k^+ - \mathbf{y},-\nabla f_{\mathbf{y}}(\mathbf{z}_k^+)\right\rangle-f(\mathbf{x}_k^+,\mathbf{y}^*)+\left\langle \mathbf{x}_k^+-\mathbf{x}^*,\nabla f_{\mathbf{x}}(\mathbf{z}_k^+)\right\rangle + f(\mathbf{x}^*,\mathbf{y}_k^+)\nonumber \\
\leq &\ \left\langle \mathbf{y}_k^+ - \mathbf{y},-\nabla f_{\mathbf{y}}(\mathbf{z}_k^+)\right\rangle + \left\langle \mathbf{x}_k^+-\mathbf{x}^*,\nabla f_{\mathbf{x}}(\mathbf{z}_k^+)\right\rangle, 
\end{align}
where the first inequality follows from the concavity of $f$ in $\mathbf{y}$, the second inequality follows from the convexity of $f$ in $\mathbf{x}$, and the last inequality follows from that $\mathbf{z}^*$ is the saddle point, which indicates $f(\mathbf{x}^*,\mathbf{y}_k^+)\leq f(\mathbf{x}^*,\mathbf{y}^*)\leq f(\mathbf{x}_k^+,\mathbf{y}^*)$.
    
Denote $\tilde{\mathbf{z}}:=(\mathbf{x}^*,\mathbf{y})$, taking the total expectation of \eqref{eq: cons_violation} and then taking the telescoping sum lead to:
\begin{align*}
&\ \mathbb{E}\left[\left\langle \mathbf{y}-\mathbf{y}^*,\sum_{k=0}^{K-1}\nabla f_{\mathbf{y}}(\mathbf{z}_k^+)\right\rangle\right]\\
\leq &\ \sum_{k=0}^{K-1}\mathbb{E}\left[\left\langle \mathbf{y}_k^+ - \mathbf{y},-\nabla f_{\mathbf{y}}(\mathbf{z}_k^+)\right\rangle + \left\langle \mathbf{x}_k^+-\mathbf{x}^*,\nabla f_{\mathbf{x}}(\mathbf{z}_k^+)\right\rangle\right] \\
=&\ \sum_{k=0}^{K-1}\mathbb{E}\left[\left\langle \mathbf{z}_k^+ - \tilde{\mathbf{z}},F(\mathbf{z}_k^+)\right\rangle \right]\\
\leq &\ \Gamma_K.
\end{align*}
Define $\mathbf{s}:=\mathbb{E}\left[\sum_{k=0}^{K-1}\nabla f_{\mathbf{y}}(\mathbf{x}_k^+,\mathbf{y}_k^+)\right]$. By the convexity of $f$ in $\mathbf{x}$, we have $\mathbf{s}\geq K\mathbb{E}\left[\nabla f_{\mathbf{y}}(\hat{\mathbf{x}}_K,\hat{\mathbf{y}}_K)\right]=K\mathbb{E}\left[\boldsymbol{\phi}(\hat{\mathbf{x}}_K)\right]$. If $\left[\mathbf{s}\right]_+=0$, we have $[\mathbb{E}\left[\boldsymbol{\phi}(\hat{\mathbf{x}}_K)\right]]_+=0$ and
the constraint violation bound 
holds. Otherwise, define $\Bar{\mathbf{y}}=\mathbf{y}^*+R\frac{\left[\mathbf{s}\right]_+}{\Vert\left[\mathbf{s}\right]_+\Vert}$. By the definition of $R$, $\bar{\mathbf{y}}$ lies within $\mathcal{Y}$. Given that $\left\langle \Bar{\mathbf{y}}-\mathbf{y}^*,\mathbb{E}\left[\sum_{k=0}^{K-1}\nabla f_{\mathbf{y}}(\mathbf{x}_k^+,\mathbf{y}_k^+)\right]\right\rangle = R \Vert\left[\mathbf{s}\right]_+\Vert \leq \Gamma_K$, we have:
\begin{equation*}
\Vert\left[\mathbb{E}\left[ \boldsymbol{\phi}(\hat{\mathbf{x}}_K)\right]\right]_+\Vert\leq \frac{\Gamma_K}{RK}.
\end{equation*}
    
As for the optimal gap of the objective function value, we apply the convexity of the objective function $\phi_0$ in $\mathbf{x}$ to get:
    \begin{align*}
        &\ \phi_0(\hat{\mathbf{x}}_K)-\phi_0^*\\
            \leq &\ \frac{1}{K}\sum_{k=0}^{K-1} \left[\phi_0(\mathbf{x}_k^+)-\phi_0^*\right]\\
            = &\ \frac{1}{K}\sum_{k=0}^{K-1} \left[f(\mathbf{x}_k^+,\mathbf{y}^*)-\left\langle \mathbf{y}^*, \nabla f_{\mathbf{y}}(\mathbf{z}_k^+)\right\rangle-\phi_0^*\right]\\
            \leq &\ \frac{1}{K}\sum_{k=0}^{K-1} \left[f(\mathbf{x}_k^+,\mathbf{y}_k^+)+\left\langle \mathbf{y}_k^+, -\nabla f_{\mathbf{y}}(\mathbf{z}_k^+)\right\rangle-\phi_0^*\right]\\
            \leq  &\ \frac{1}{K}\sum_{k=0}^{K-1} \left[\left\langle \mathbf{x}_k^+-\mathbf{x}^*,\nabla f_{\mathbf{x}}(\mathbf{z}_k^+)\right\rangle + f(\mathbf{x}^*,\mathbf{y}_k^+)\right]+\frac{1}{K}\sum_{k=0}^{K-1} \left[\left\langle \mathbf{y}_k^+, -\nabla f_{\mathbf{y}}(\mathbf{z}_k^+)\right\rangle-f(\mathbf{x}^*,\mathbf{y}^*)\right]\\
            \leq &\ \frac{1}{K}\sum_{k=0}^{K-1} \left[\left\langle \mathbf{x}_k^+-\mathbf{x}^*,\nabla f_{\mathbf{x}}(\mathbf{z}_k^+)\right\rangle \right]+\frac{1}{K}\sum_{k=0}^{K-1} \left[
            \left\langle \mathbf{y}_k^+, -\nabla f_{\mathbf{y}}(\mathbf{z}_k^+)\right\rangle\right],
    \end{align*}
where the first inequality follows from Jensen's inequality; the second inequality follows from the concavity of $f$ in $\mathbf{y}$; the third inequality follows from the convexity of $f$ in $\mathbf{x}$; and the last one is derived from $f(\mathbf{x}^*,\mathbf{y}_k^+)\leq f(\mathbf{x}^*,\mathbf{y}^*)$.
Let $\check{\mathbf{z}}:=(\mathbf{x}^*,\mathbf{0})\in\mathcal{Z}$, then we have:

$$\sum_{k=0}^{K-1}\mathbb{E}\left[\left\langle F(\mathbf{z}_k^+),\mathbf{z}_k^+-\check{\mathbf{z}}\right\rangle\right]
    \leq \Gamma_K.$$
Consequently, the optimal gap of $\phi_0(\hat{\mathbf{x}}_K)$ can be bounded by $\mathbb{E}\left[\phi_0(\hat{\mathbf{x}}_K)\right]-\phi_0^*\leq \frac{\Gamma_K}{K}$.
\end{pf*}
Substituting \eqref{eq: bound_ZOEG_telescope} into Lemma~\ref{lem: optimality gap} finishes the proof.
\end{pf*}

It is shown in Proposition~\ref{prop: ZOEG vio obj} that the constraint violation and the optimality gap can both converge at the same convergence rate with the duality gap. For the bounds in Theorem \ref{thm: ZOEG_duality gap} and Proposition \ref{prop: ZOEG vio obj}, we can apply small step sizes to eliminate the fixed error terms.
As a corollary, the iteration complexity required to obtain an $\epsilon$-saddle point for problem \eqref{eq: min-max} and an $\epsilon$-optimal solution for problem \eqref{eq: constrained_problem} is provided as follows.

\begin{corollary}\label{cor: ZOEG_complexity}
Suppose the assumptions and conditions in Theorem \ref{thm: ZOEG_duality gap} hold. Set $\eta=\frac{\tilde{D}}{6\sqrt{6dK}G}$ and $K\geq \max\left\{ \frac{864 d G^2 \tilde{D}^2}{\epsilon^2},\left(\frac{(6d^2+8)^2L^4M_2^2\tilde{D}^2)}{6dG^2\epsilon^2}\right)^{\frac{1}{3}},\frac{4LM_1\tilde{D}}{\epsilon}\right\}$
Then, we have:
\begin{equation*}
\mathbb{E}\left[\Delta(\hat{\mathbf{z}}_K)\right]\leq \epsilon,
\end{equation*} 
\begin{equation*}
\Vert\left[\mathbb{E}\left[ \boldsymbol{\phi}(\hat{\mathbf{x}}_K)\right]\right]_+\Vert\leq \frac{\epsilon}{R},
    \end{equation*}
    \begin{equation*}
\mathbb{E}\left[\phi_0(\hat{\mathbf{x}}_K)\right]-\phi_0^*\leq \epsilon.
    \end{equation*} 
\end{corollary}
\begin{pf*}[Proof of Corollary \ref{cor: ZOEG_complexity}]
By the results in Theorem \ref{thm: ZOEG_duality gap}, substituting $\eta=\frac{\tilde{D}}{6\sqrt{6dK}G}$ into $\frac{\frac{\tilde{D}^2}{2\eta }+C_1 \eta + C_2}{K}+108\eta d G^2$ gives:
\begin{align*}
&\ \mathbb{E}\left[\Delta(\hat{\mathbf{z}}_K)\right]\\
\leq &\ \frac{\tilde{D}^2}{2\eta K}+108\eta d G^2+\frac{(9d^2+12)\eta L^2 M_2+L\tilde{D}M_1}{K}\\
=&\ 6\sqrt{\frac{6d}{K}}\tilde{D}G+\frac{(3d^2+4)L^2M_2\tilde{D}}{2\sqrt{6dK^3}G}+\frac{LM_1\tilde{D}}{K}
            \\
\leq &\ \frac{\epsilon}{2}+\frac{\epsilon}{4}+\frac{\epsilon}{4}\\
=&\ \epsilon,
\end{align*}
where the first equality follows from $\eta=\frac{\tilde{D}}{6\sqrt{6dK}G}$ and the last inequality follows from the lower bound of $K$. The remaining two bounds follow a similar proof.
\end{pf*}

The above corollary indicates that the iteration complexity to achieve an $\epsilon$-approximate solution for problems \eqref{eq: constrained_problem} and \eqref{eq: min-max} is $\mathcal{O}\left(\frac{d}{\epsilon^2}\right)$, which is also the oracle complexity. The dependence on $\epsilon$ matches the best-known result of $\mathcal{O}(\frac{d^2}{\epsilon^2})$ in \cite{nguyen2023stochastic} for solving problem \eqref{eq: constrained_problem} and another result of $\mathcal{O}(\frac{d^4}{\epsilon^2})$ in \cite{maheshwari2022zeroth} for solving problem \eqref{eq: min-max}. Moreover, our results achieve a better and the best-known dependence on $d$.
\begin{remark}
While our analysis focuses on the deterministic setting, ZOEG and the main theoretical results are also applicable to the stochastic case, as discussed in \cite{nguyen2023stochastic}. Stochasticity introduces an additional variance term in the gradient estimation bound (see \eqref{eq: var_2p_grad}), but our complexity results remain unchanged.
\end{remark}
Recall that the oracle complexity of the first-order EG is $\mathcal{O}(\frac{1}{\epsilon})$ \cite{nemirovski2004prox}, which means ZOEG is strictly inferior to the first-order method. Therefore, it is very natural to ask whether we can improve ZOEG to bridge the gap. In the next section, we provide a positive answer to this question by leveraging the $2d$-point coordinate gradient estimator.

\section{Zeroth-Order Coordinate Extra-Gradient}\label{sec: ZOCEG}
In this section, we will present another algorithm called zeroth-order coordinate extra-gradient (ZOCEG) and provide its convergence guarantees. 
\subsection{Algorithm Design}
According to the analysis in Section~\ref{sec: ZOEG}, ZOEG suffers from slow convergence. The main factor that limits the performance of ZOEG is the large variance of $2$-point UniGE \eqref{eq: 2p_grad_est}. By Lemma \ref{lem: coor_grad_estimator}, the coordinate-wise partial gradient estimator~\eqref{eq: coor_grad_est} enjoys controllable bias and variance by using small smoothing radii.
To leverage such merits, we integrate $2d$-point CooGE into the EG method and develop the ZOCEG algorithm. The gradient estimators of $\nabla_{\mathbf{x}}f(\mathbf{x}_k,\mathbf{y}_k)$ and $\nabla_{\mathbf{y}}f(\mathbf{x}_k,\mathbf{y}_k)$ are given by:
\begin{align*}
\hat{\mathbf{g}}^{x}_k=\sum_{i=1}^{d_x}\frac{f(\mathbf{x}_k+r_k\mathbf{e}_i,\mathbf{y}_k)-f(\mathbf{x}_k,\mathbf{y}_k)}{r_k}\mathbf{e}_i,\\
\hat{\mathbf{g}}^{y}_k=\sum_{j=1}^{d_y}\frac{f(\mathbf{x}_k,\mathbf{y}_k+r_k\mathbf{e}_j)-f(\mathbf{x}_k,\mathbf{y}_k)}{r_k}\mathbf{e}_j.
\end{align*}
With $\mathbf{x}_k$ and $\mathbf{y}_k$ updated by the estimated gradients, we get $(\mathbf{x}_k^+,\mathbf{y}_k^+)$, and the gradient of $f(\mathbf{x}_k^+,\mathbf{y}_k^+)$ can be estimated similarly as follows:
\begin{align*}
\hat{\mathbf{g}}^{x,+}_k=\sum_{i=1}^{d_x}\frac{f(\mathbf{x}^+_k+r_k\mathbf{e}_i,\mathbf{y}^+_k)-f(\mathbf{x}^+_k,\mathbf{y}^+_k)}{r_k}\mathbf{e}_i,\\
\hat{\mathbf{g}}^{y,+}_k=\sum_{j=1}^{d_y}\frac{f(\mathbf{x}^+_k,\mathbf{y}^+_k+r_k\mathbf{e}_j)-f(\mathbf{x}^+_k,\mathbf{y}^+_k)}{r_k}\mathbf{e}_j.
\end{align*}
The detailed procedure is presented in Algorithm~\ref{alg: ZOCEG}.
\begin{algorithm}[htbp]
	\caption{Zeroth-Order Coordinate Extra-Gradient (ZOCEG)}
\begin{algorithmic}\label{alg: ZOCEG}
	\STATE \textbf{Inputs:} Number of iterations $K$, step size $\eta$, smoothing radii $\{r_k\}$, and initial point $\mathbf{z}_0$.
	\FOR{$k\xleftarrow{}0$ to $K-1$} 
        \STATE Calculate $\mathbf{z}_k^+$ by:
        \begin{equation*}\label{eq: ZOCEG_plus_x}
            \mathbf{x}^+_k=\mathcal{P}_{\mathcal{X}}\left [\mathbf{x}_k-\eta\hat{\mathbf{g}}_k^{x}\right],
        \end{equation*}
        \begin{equation*}\label{eq: ZOCEG_plus_y}
            \mathbf{y}^+_k=\mathcal{P}_{\mathcal{Y}}\left [\mathbf{y}_k+\eta\hat{\mathbf{g}}_k^{y}\right].
        \end{equation*}
        Based on $\mathbf{z}_k^+$, update $\mathbf{z}_k$ by:
        \begin{equation*}\label{eq: ZOCEG_next_x}
            \mathbf{x}_{k+1}=\mathcal{P}_{\mathcal{X}}\left [\mathbf{x}_k-\eta\hat{\mathbf{g}}_k^{x,+}\right],
        \end{equation*}
        \begin{equation*}\label{eq: ZOCEG_next_y}
            \mathbf{y}_{k+1}=\mathcal{P}_{\mathcal{Y}}\left [\mathbf{y}_k+\eta\hat{\mathbf{g}}_k^{y,+}\right].
        \end{equation*}
	\ENDFOR
    \STATE \textbf{Outputs:} $\hat{\mathbf{z}}_K=\frac{1}{K}\sum_{k=0}^{K-1} \mathbf{z}^+_k .$
\end{algorithmic}
\end{algorithm}
\subsection{Finite-Sample Analysis}
First, we focus on the performance of ZOCEG when solving problem \eqref{eq: min-max}. The duality gap of the averaged iterate from ZOCEG is bounded in the following theorem.
\begin{theorem}\label{thm: ZOCEG duality gap}
  Suppose Assumptions \ref{ass: convexity}-\ref{ass: bounded} hold. 
 The sequence $\{\mathbf{z}_k^+\}$ is generated by the ZOCEG Algorithm. Set the smoothing radii to satisfy 
$\sum_{k=1}^{\infty}r_k\leq\frac{M_3}{\sqrt{d_x}+\sqrt{d_y}}<\infty$ for some $M_3>0$. Let the step sizes satisfy $0<\eta L\leq \frac{1}{2}$, then for any $K>1$, we have:
    \begin{equation*}
\Delta(\hat{\mathbf{z}}_K)\leq \frac{1}{2K}\left(\frac{\tilde{D}^2}{\eta}+3 LM_3\tilde{D}\right).
    \end{equation*}
\end{theorem}
\begin{pf*}[Proof of Theorem \ref{thm: ZOCEG duality gap}]
Following a similar idea to analyzing ZOEG, we start from providing the bound of $\left\langle\eta F(\mathbf{z}_k^+),\mathbf{z}_k^+-\mathbf{z}\right\rangle$ in the following lemma, the proof of which is deferred to Appendix~\ref{pf: bound ZOCEG term}. 
\begin{lemma}\label{lem: bound ZOCEG term}
    Suppose Assumptions~\ref{ass: convexity}-\ref{ass: bounded} hold and the step sizes satisfy $0<\eta L\leq \frac{1}{2}$. Given the iterates from ZOCEG, we have for any $\mathbf{z}\in\mathcal{Z}$:
\begin{equation*} 
    \left\langle\eta F(\mathbf{z}_k^+),\mathbf{z}_k^+-\mathbf{z}\right\rangle
    \leq H_k(\mathbf{z})-H_{k+1}(\mathbf{z})+\frac{3}{2}\eta Lr_k(\sqrt{d_x}+\sqrt{d_y})\tilde{D}.
\end{equation*}
\end{lemma}
Based on Lemma~\ref{lem: bound ZOCEG term}, taking the telescoping sum of $\left\langle F(\mathbf{z}_k^+),\mathbf{z}_k^+-\mathbf{z}\right\rangle$ gives the upper bound as follows:
\begin{align}\label{eq: bound_ZOCEG_telescope}
    &\ \sum_{k=0}^{K-1}\left\langle F(\mathbf{z}_k^+),\mathbf{z}_k^+-\mathbf{z}\right\rangle\nonumber\\
    \leq&\ \frac{1}{\eta}H_{0}(\mathbf{z})-\frac{1}{\eta}H_{K}(\mathbf{z})+\sum_{k=0}^{K-1}\frac{3}{2} Lr_k(\sqrt{d_x}+\sqrt{d_y})\tilde{D}\nonumber\\
    \leq&\ \frac{\tilde{D}^2}{2\eta}+\frac{3}{2} L
\frac{M_3}{\sqrt{d_x}+\sqrt{d_y}}\cdot(\sqrt{d_x}+\sqrt{d_y})\tilde{D}\nonumber\\
    \leq& \frac{\tilde{D}^2}{2\eta}+\frac{3}{2} LM_3\tilde{D},
\end{align}
where the last inequality follows from the definition of $H_k(\mathbf{z})$ and the boundedness of $\mathcal{Z}$.
Substituting \eqref{eq: bound_ZOCEG_telescope} into Lemma~\ref{lem: duality gap}, which still holds without the expectation, finishes the proof.
\end{pf*}

The result in Theorem~\ref{thm: ZOCEG duality gap} shows that the generated iterates from ZOCEG can converge to a saddle point at the convergence rate of $\mathcal{O}\left(\frac{1}{K}\right)$. Following the same idea in the analysis of ZOEG, combining \eqref{eq: bound_ZOCEG_telescope} with the deterministic counterpart of Lemma~\ref{lem: optimality gap} directly establishes the following proposition, which provides finite-iteration bounds on the optimality gap and constraint violation of ZOCEG for solving the original constrained problem \eqref{eq: constrained_problem}. Therefore, we omit its proof here.
\begin{proposition}\label{prop: ZOCEG vio obj}
Suppose the Assumptions and conditions in Theorem \ref{thm: ZOCEG duality gap} hold. Then for any $K>1$, we have:
\begin{align*}
\Vert\left[\boldsymbol{\phi}(\hat{\mathbf{x}}_K)\right]_+\Vert &\leq \frac{\tilde{D}^2}{2RK\eta}+\frac{3LM_3\tilde{D}}{2RK},\\
\phi_0(\hat{\mathbf{x}}_K)-\phi_0^* &\leq \frac{\tilde{D}^2}{2K\eta}+\frac{3LM_3\tilde{D}}{2K}.
\end{align*}
\end{proposition}

Based on Theorem~\ref{thm: ZOCEG duality gap} and Proposition~\ref{prop: ZOCEG vio obj}, we can further establish the oracle complexity of ZOCEG when solving problems \eqref{eq: constrained_problem} and \eqref{eq: min-max}.
\begin{corollary}\label{cor: ZOCEG complexity}
Suppose the Assumptions and conditions in Theorem \ref{thm: ZOCEG duality gap} hold. Then for any $ K\geq \max\left\{ \frac{\tilde{D}^2}{\eta\epsilon},\frac{3LM_3\tilde{D}}{\epsilon}\right\}$, we have
\begin{equation*}
    \Delta(\hat{\mathbf{z}}_K)\leq \epsilon,
\end{equation*}
\begin{equation*}
    \Vert\left[ \boldsymbol{\phi}(\hat{\mathbf{x}}_K)\right]^+\Vert \leq \frac{\epsilon}{R},
\end{equation*}
\begin{equation*}
    \mathbb{E}\left[\phi_0(\hat{\mathbf{x}}_K)\right]-\phi_0^* \leq \epsilon.
\end{equation*}
\end{corollary}

The proof of Corollary \ref{cor: ZOCEG complexity} follows the same idea as Corollary~\ref{cor: ZOEG_complexity} and is thus omitted. Corollary~\ref{cor: ZOCEG complexity} shows that the oracle complexity for ZOCEG to derive an $\epsilon$-saddle point for problem \eqref{eq: min-max} and an $\epsilon$-optimal solution for problem \eqref{eq: constrained_problem} is $\mathcal{O}\left(\frac{d}{\epsilon}\right)$, which differs from first-order EG by a factor of $d$ \cite{nemirovski2004prox}.
This is the best-known complexity bound for zeroth-order algorithms when solving problem \eqref{eq: constrained_problem}. 

\begin{remark}
It is worth noting that the complexity bounds proportional to $\epsilon^{-1}$ are also derived in \cite{he2022zeroth} and \cite{jin2023zeroth} for deterministic constrained ZO algorithms. However, both studies rely on a decomposable and explicit feasible set, which cannot apply to problem \eqref{eq: constrained_problem} with black-box constraints.
\end{remark}
\begin{remark}
Compared to ZOEG, the oracle complexity is improved from $\mathcal{O}\left(\frac{d}{\epsilon^2}\right)$ to $\mathcal{O}\left(\frac{d}{\epsilon}\right)$. The key factor driving the improvement is that the bias and variance of the $2d$-point CooGE can be controlled to be sufficiently small, which makes ZOCEG resemble the behavior of first-order EG. However, ZOCEG and its theoretical results do not apply to the stochastic setting, where extra variance is introduced to gradient estimation due to the system's stochasticity.
\end{remark}
\begin{remark}
We consider the convex problems throughout this work. As for nonconvex problems, we can combine our algorithms with the PP method in \cite{nguyen2023stochastic} and \cite{boob2023stochastic}, a meta-algorithm, to design a double-loop algorithm. Using the theoretical results in \cite{boob2023stochastic}, it is not hard to derive the oracle complexities of $\mathcal{O}\left(\frac{d}{\epsilon^3}\right)$ and $\mathcal{O}\left(\frac{d}{\epsilon^2}\right)$ for ZOEG and ZOCEG-based PP algorithms to get an $\epsilon$-critical KKT point in nonconvex problems. The result of $\mathcal{O}\left(\frac{d}{\epsilon^2}\right)$ also matches recent studies for nonconvex ZO problems \cite{guo2025safe,xu2024derivative}. When directly applying our algorithms to nonconvex problems, we leave the complexity analysis for future work and only conduct empirical validation in this study.
\end{remark}
\subsection{Zeroth-Order Block Coordinate Extra-Gradient}
While ZOCEG has excellent oracle complexity guarantees, each gradient estimation requires $\mathcal{O}(d)$ function values, which is inefficient for high-dimensional problems.
To improve computational efficiency for each single step, we may further consider only estimating partial gradients along a block of coordinates rather than the full gradient and updating the block of coordinates accordingly. Let $\tau_x,\tau_y$ denote the number of coordinates updated for the primal and dual variables at each iteration, where $1\leq \tau_x\leq d_x$, and $1\leq \tau_y\leq d_y$. At the $k$th iteration, the coordinate index subsets $\mathcal{I}_k^x\subseteq\{1,2,\cdots, d_x\}$ and $\mathcal{I}_k^y\subseteq\{1,2,\cdots, d_y\}$, with $|\mathcal{I}_k^x|=\tau_x$, $|\mathcal{I}_k^y|=\tau_y$, are sampled uniformly without replacement.
The partial gradient vectors 
$\left\{\nabla_{\mathbf{x},i}f(\mathbf{x}_k,\mathbf{y}_k)\right\}_{i\in\mathcal{I}_k^x}$  and $\left\{\nabla_{\mathbf{y},j}f(\mathbf{x}_k,\mathbf{y}_k)\right\}_{j\in\mathcal{I}_k^y}$  corresponding to selected coordinates can then be estimated as follows:
\begin{align*}
\tilde{\mathbf{g}}^{x}_k=\sum_{i\in\mathcal{I}_k^x}\frac{f(\mathbf{\mathbf{x}}_k+r_k\mathbf{e}_i,\mathbf{y}_k)-f(\mathbf{\mathbf{x}}_k,\mathbf{y}_k)}{r_k}\mathbf{e}_i,\\
\tilde{\mathbf{g}}^{y}_k=\sum_{j\in\mathcal{I}_k^y}\frac{f(\mathbf{\mathbf{x}}_k,\mathbf{y}_k+r_k\mathbf{e}_j)-f(\mathbf{\mathbf{x}}_k,\mathbf{y}_k)}{r_k}\mathbf{e}_j.
\end{align*}

Using the 2$\tau_x$-point CooGE and 2$\tau_y$-point CooGE, only the entries of $\mathbf{x}_k$ indexed by $\mathcal{I}_k^x$ and the entries of $\mathbf{y}_k$ indexed by $\mathcal{I}_k^y$ are updated, yielding the intermediate variables $(\mathbf{x}_k^+,\mathbf{y}_k^+)$.
Similarly, another pair of coordinate index subsets $\mathcal{J}_k^x$ and $\mathcal{J}_k^y $ are independently sampled, and the corresponding partial 
gradient vectors at $(\mathbf{x}_k^+,\mathbf{y}_k^+)$ are estimated as:
\begin{align*}
\tilde{\mathbf{g}}^{x,+}_k=\sum_{p\in\mathcal{J}_k^x}\frac{f(\mathbf{\mathbf{x}}^+_k+r_k\mathbf{e}_p,\mathbf{y}_k^+)-f(\mathbf{\mathbf{x}}_k^+,\mathbf{y}_k^+)}{r_k}\mathbf{e}_p,\\
\tilde{\mathbf{g}}^{y,+}_k=\sum_{q\in\mathcal{J}_k^y}\frac{f(\mathbf{\mathbf{x}}^+_k,\mathbf{y}_k^++r_k\mathbf{e}_q)-f(\mathbf{\mathbf{x}}^+_k,\mathbf{y}^+_k)}{r_k}\mathbf{e}_q,
\end{align*}
which are then used to compute the new iterates $\left(\mathbf{x}_{k+1},\mathbf{y}_{k+1}\right)$.

The detailed algorithm is presented in Algorithm~\ref{alg: ZOBCEG}. 
Note that the blocks for gradient estimation of $(\mathbf{x}_k,\mathbf{y}_k)$ and $(\mathbf{x}_k^+,\mathbf{y}_k^+)$ are sampled independently.  The convergence analysis of ZOBCEG is provided in Theorem~\ref{thm: ZOBCEG convergence}, the proof of which is deferred to Appendix~\ref{apx: pf_ZOBCEG}.

\begin{theorem}\label{thm: ZOBCEG convergence}
  Suppose Assumptions \ref{ass: convexity}-\ref{ass: bounded} hold. 
 The sequence $\{\mathbf{z}_k^+\}$ is generated by the ZOBCEG Algorithm. Set the smoothing radii to satisfy 
$\sum_{k=1}^{\infty}r_k=M_1<\infty$ and $\sum_{k=1}^{\infty}r_k^2=M_2<\infty$
, then for any $K>1$, we have:
\begin{align*}    \mathbb{E}\left[\Delta(\hat{\mathbf{z}}_K)\right]\leq\left(\frac{d_x}{\tau_x}+\frac{d_y}{\tau_y}\right)\frac{\tilde{D}^2}{2\eta K}+\frac{LM_1(\sqrt{d_x}+\sqrt{d_y})\tilde{D}}{2K}+\frac{3(d_x+d_y)}{2K}\eta L^2 M_2+24\eta G^2.
\end{align*}

\end{theorem}
Based on the result in Theorem \ref{thm: ZOBCEG convergence}, ZOBCEG can achieve the convergence to a saddle point at the rate of $\mathcal{O}\left(\frac{\tilde{d}}{\eta K}+\frac{\eta d}{K}+\eta \right)$, where $\tilde{d}=\max\left\{\frac{d_x}{\tau_x},\frac{d_y}{\tau_y}\right\}$. Similarly to ZOEG, we can apply a small step size of $\eta$, i.e., $\eta=\mathcal{O}\left(\sqrt{\frac{\tilde{d}}{K}}\right)$, to achieve the convergence rate of $\mathcal{O}\left(\frac{\sqrt{\tilde{d}}}{\sqrt{K}}\right)$. Furthermore, the oracle complexity of ZOBCEG is $\mathcal{O}\left(\frac{d}{\epsilon^2}\right)$ to derive an $\epsilon$-saddle point for problem \eqref{eq: min-max} and an $\epsilon$-optimal solution for problem \eqref{eq: constrained_problem}. 

Compared with ZOCEG, the theoretical guarantees of ZOBCEG degrade due to the random selection of a block of coordinates for updates, which introduces further variance for gradient estimation. Therefore, we can only derive similar theoretical guarantees with ZOEG with improved dependence on the dimension. When $\tau_x=d_x$, and $\tau_y=d_y$, ZOBCEG estimates full gradients and recovers the ZOCEG algorithm. However, our subsequent numerical results show that ZOBCEG can typically achieve comparable oracle complexities to ZOCEG even when $\tau_x<d_x$ and $\tau_y<d_y$.
\begin{algorithm}[htbp]
	\caption{Zeroth-Order Block Coordinate Extra-Gradient (ZOBCEG)}
\begin{algorithmic}\label{alg: ZOBCEG}
	\STATE \textbf{Inputs:} Number of iterations $K$, step size $\eta$, smoothing radii $\{r_k\}$, and initial point $\mathbf{z}_0$.
	\FOR{$k\xleftarrow{}0$ to $K-1$} 
        \STATE Randomly sample $\mathcal{I}_k^x\subseteq\{1,2,\cdots, d_x\}$,  $\mathcal{I}_k^x\subseteq\{1,2,\cdots, d_x\}$, with $|\mathcal{I}_k^x|=\tau_x$, $|\mathcal{I}_k^y|=\tau_y$, and calculate $\mathbf{z}_k^+$ by:
        \begin{equation*}\label{eq: ZORCEG_plus_x}
\mathbf{x}^+_k=\mathcal{P}_{\mathcal{X}}\left [\mathbf{x}_k-\eta\tilde{\mathbf{g}}_k^{x}\right],
        \end{equation*}
        \begin{equation*}\label{eq: ZORCEG_plus_y}
\mathbf{y}^+_k=\mathcal{P}_{\mathcal{Y}}\left [\mathbf{y}_k+\eta\tilde{\mathbf{g}}_k^{y}\right].
        \end{equation*}
        Based on $\mathbf{z}_k^+$, randomly sample $\mathcal{J}_k^x\subseteq\{1,2,\cdots, d_x\}$,  $\mathcal{J}_k^x\subseteq\{1,2,\cdots, d_x\}$, with $|\mathcal{J}_k^x|=\tau_x$, $|\mathcal{J}_k^y|=\tau_y$, and update $\mathbf{z}_k$ by:
        \begin{equation*}\label{eq: ZORCEG_next_x}
\mathbf{x}_{k+1}=\mathcal{P}_{\mathcal{X}}\left [\mathbf{x}_k-\eta\tilde{\mathbf{g}}_k^{x,+}\right],
        \end{equation*}
        \begin{equation*}\label{eq: ZORCEG_next_y}
\mathbf{y}_{k+1}=\mathcal{P}_{\mathcal{Y}}\left [\mathbf{y}_k+\eta\tilde{\mathbf{g}}_k^{y,+}\right].
        \end{equation*}
	\ENDFOR
    \STATE \textbf{Outputs:} $\hat{\mathbf{z}}_K=\frac{1}{K}\sum_{k=0}^{K-1} \mathbf{z}^+_k .$
\end{algorithmic}
\end{algorithm}

\section{Numerical Simulations}
In this section, we conduct numerical experiments to validate the performance of the proposed algorithms on the demand-side load tracking problem in power systems in both convex and nonconvex test cases, as in \cite{jin2023zeroth}. An aggregator tries to coordinate the loads of multiple users to curtail the total power load $p_c(\mathbf{x})$ and track the target load level while minimizing the total cost $\phi_0(\mathbf{x})$. The problem is formulated as follows:
\begin{subequations}
    \begin{align*}
        &\min_{\mathbf{x}\in\mathcal{X}}\  \phi_0(\mathbf{x}), \\
        &s.t.\ \ p_c(\mathbf{x})\leq D,
    \end{align*}
\end{subequations}
where $p_c(\mathbf{x})$ is determined by the dynamics of the distribution system and seen as a black box. 

\subsection{Convex Case}
In the convex case, we apply a simplified linear model of $p_c(x):\;p_c(\mathbf{x})=\sum_{i=1}^{100}(1+\gamma_i)(u_i-x_i)
$, where $\mathbf{x}$ denotes the load adjustment vector. We consider $\phi_0(\mathbf{x})=\sum_{i=1}^{100} h_i(x_i)$ and $ h_i(x_i)$ is the local cost function of user $i$, which is modeled as $ h_i(x_i)=a_ix_i^2+b_ix_i$. The feasible set is defined as $\mathcal{X}=\left\{\mathbf{x}\in\mathbb{R}^{100}\mid 0\leq x_i\leq u_{i},i=1,...,100\right\}$.

In addition to testing our proposed algorithms, we also compare them with two other algorithms, stochastic zeroth-order constraint extrapolation (SZO-ConEX) \cite{nguyen2023stochastic} and zeroth-order gradient descent ascent (ZOGDA) \cite{liu2020min}, applicable to solving problem \eqref{eq: constrained_problem}.
Both of them employ the two-point gradient estimator. We run the tests with constant step size (CS) and diminishing step size (DS), respectively, where the DS is set as $\eta_k:=\eta_0/\sqrt{k+1}$. Moreover, we evaluate the performance of ZOBCEG under different block sizes (BS) $\tau_x$. Note that setting $BS=100$ recovers the ZOCEG algorithm, which adopts full gradient estimation. 

\begin{figure}[htb]
    \centering
    \subfloat[]{  
\includegraphics[scale=0.4]{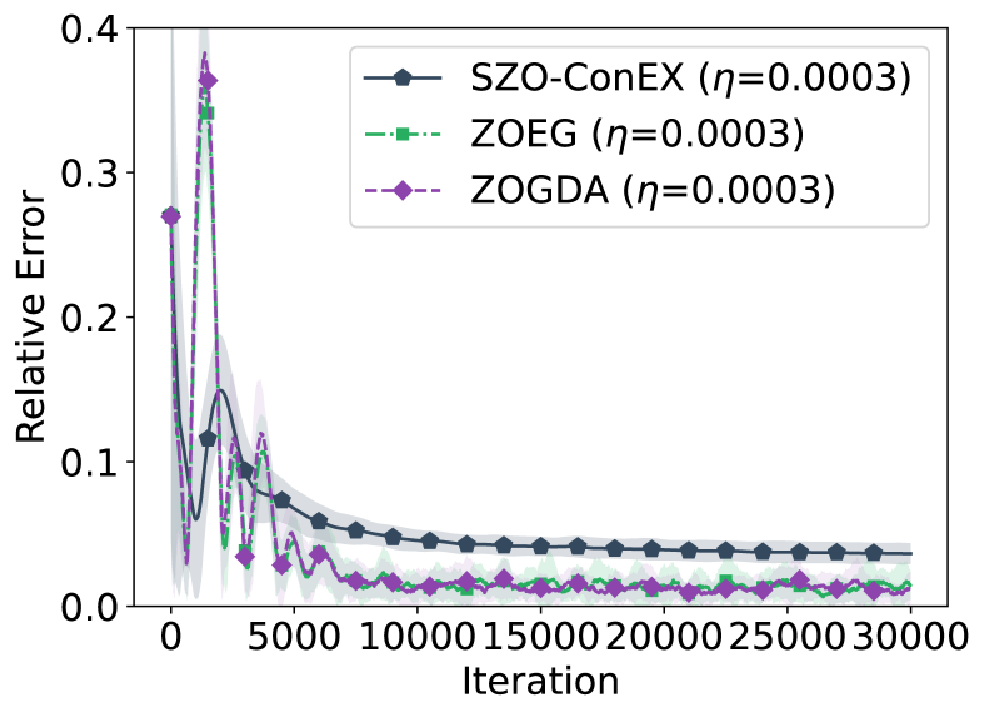}
          \label{figure convex_RE_CS}
    }   
    \hfil
    \subfloat[]{      \includegraphics[scale=0.4]{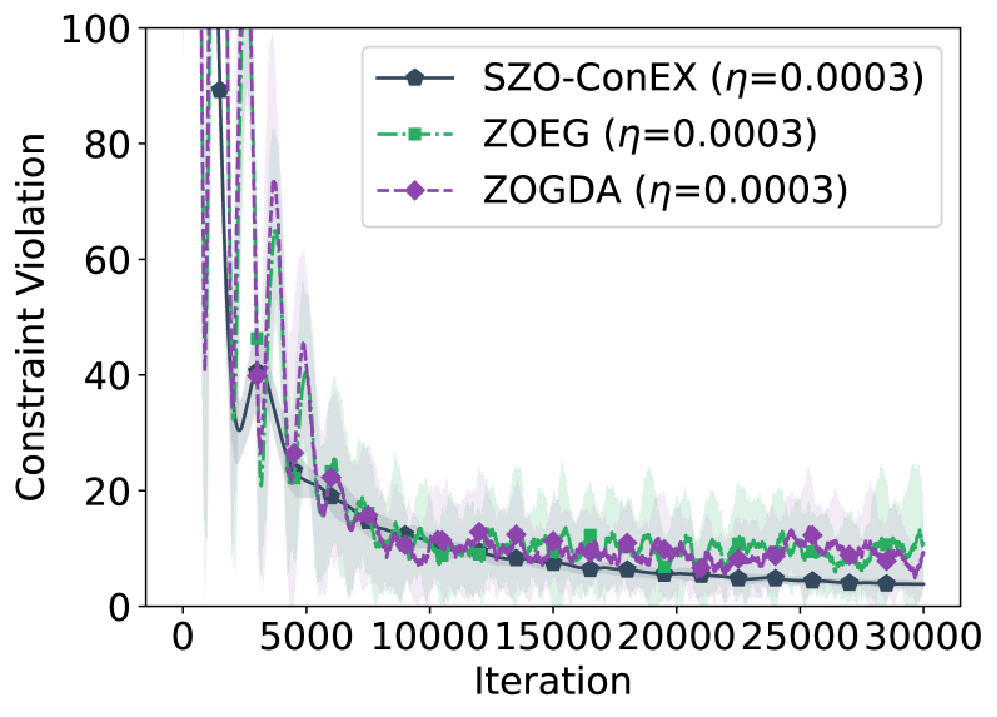}
          \label{figure convex_vio_CS}
    }
    \par\smallskip
    \subfloat[]{   \includegraphics[scale=0.4]{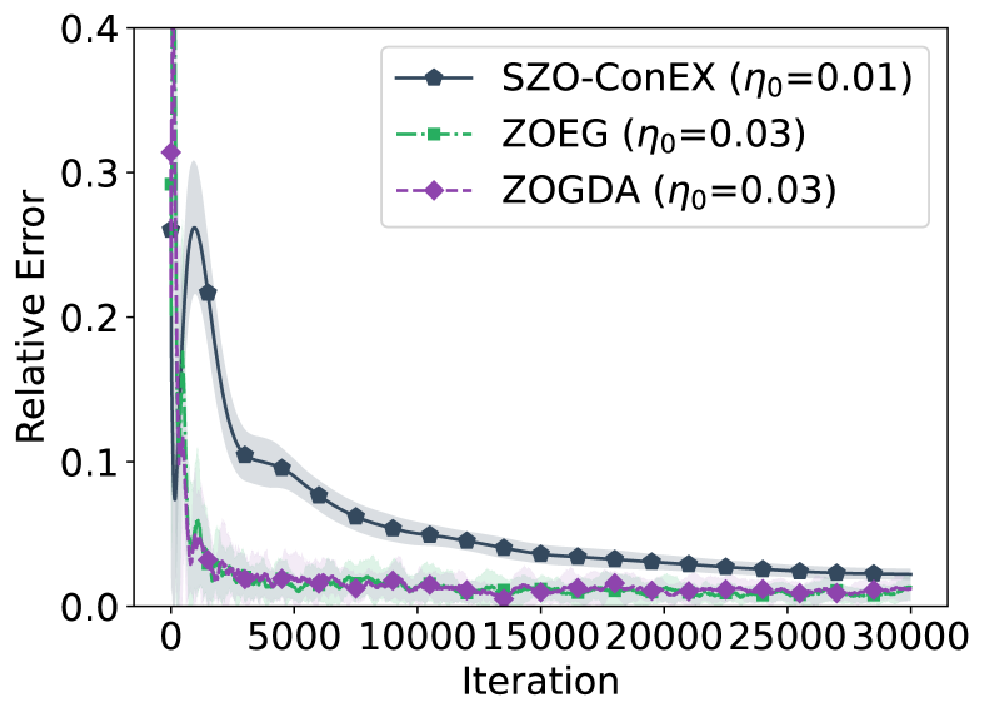}
          \label{figure convex_RE_DS}
	}
    \hfil
    \subfloat[]{   \includegraphics[scale=0.4]{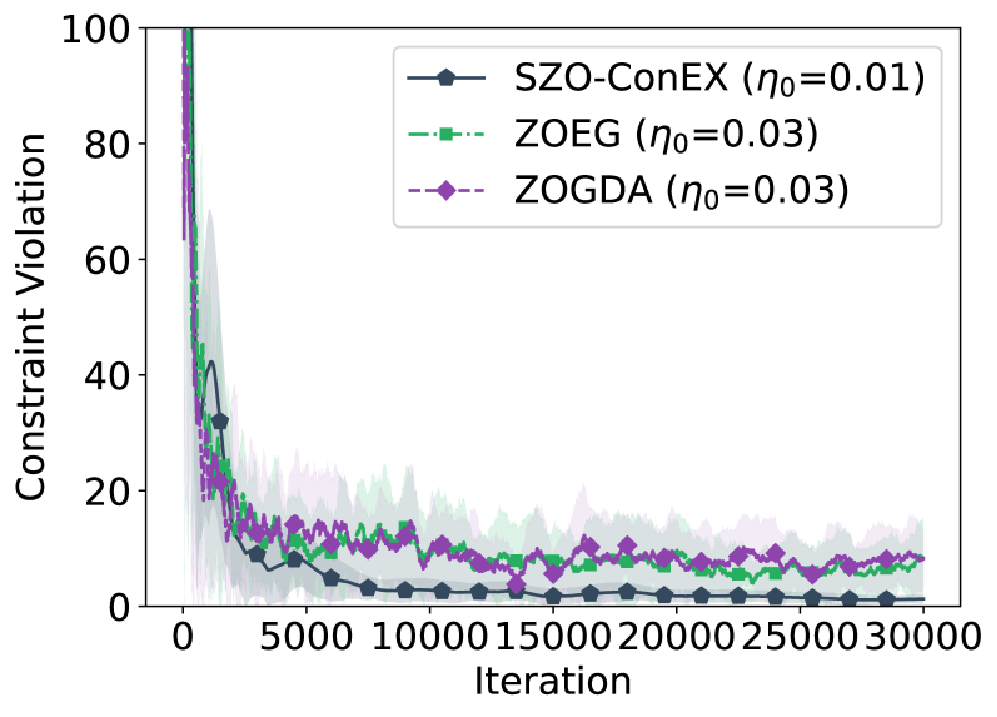}
          \label{figure convex_vio_DS}
	}
    \caption{Performance comparison in the convex case: (a) average relative error under CS, (b) average constraint violation under CS, (c) average relative error under DS, (d) average constraint violation under DS. }
    \label{figure convex performance}
\end{figure}

\begin{figure}[htb]
    \centering
    \subfloat[]{  
\includegraphics[scale=0.4]{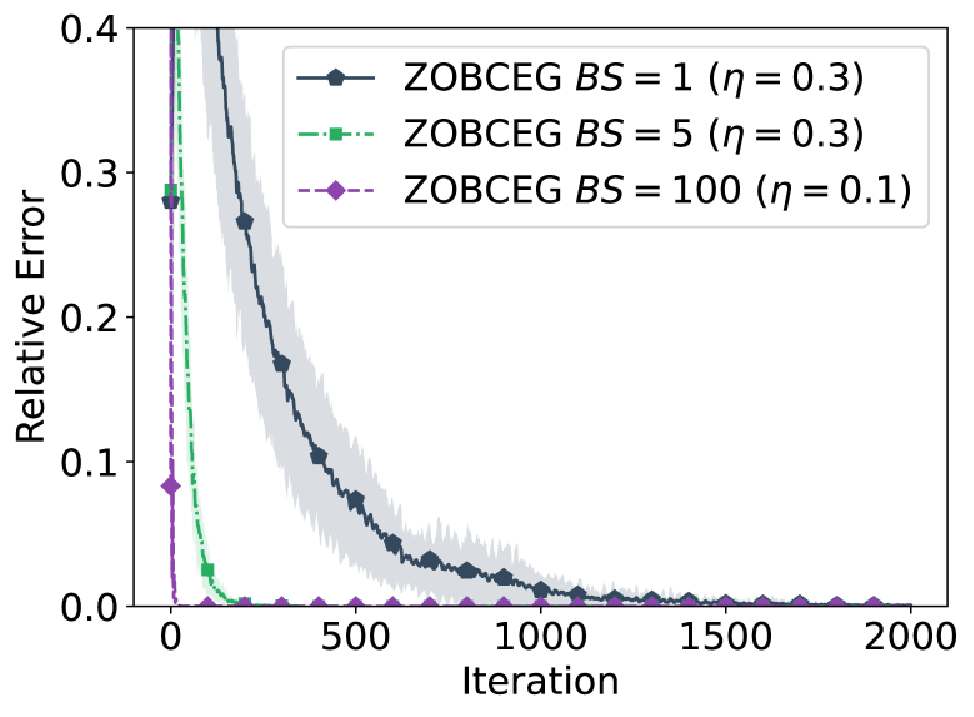}


          \label{figure convex_re_ZOBCEG}
    }   
    \hfil
    \subfloat[]{      \includegraphics[scale=0.4]{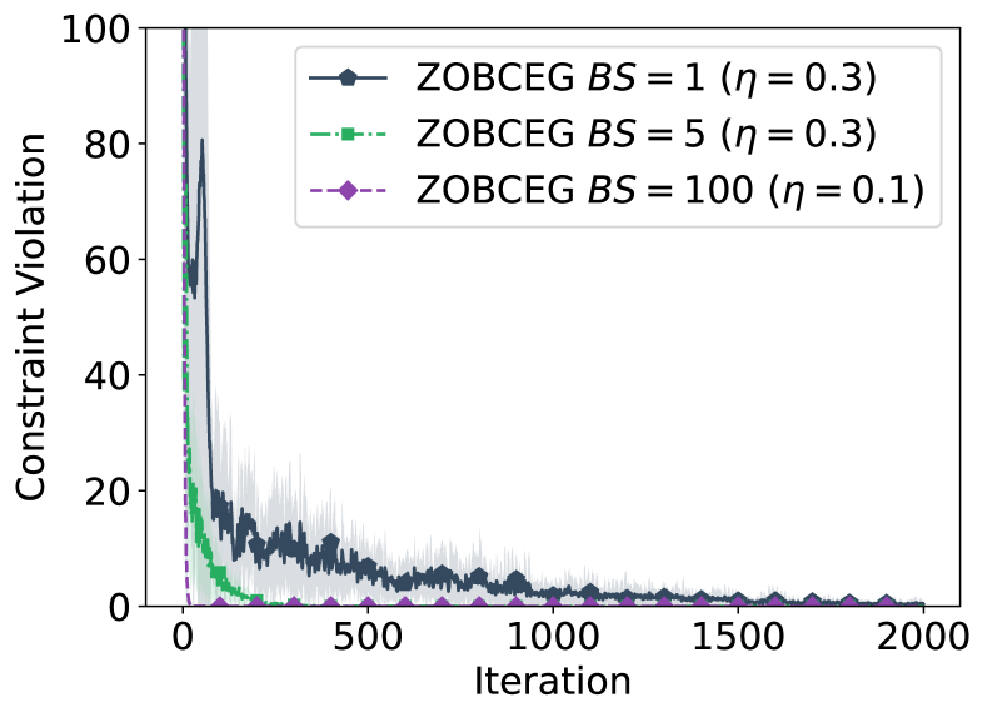}}
          \label{figure convex_vio_ZOBCEG}
    
    \caption{Performance of ZOBCEG in the convex case: (a) average relative error, (b) average constraint violation. }
    \label{figure convex ZOCEG}
\end{figure}

Regarding the experimental setup, the coefficients $a_i$ and $b_i$ in the local cost functions are sampled from the uniform distribution $a_i\sim U(0.5,1.5)$ and $b_i\sim U(0,5)$, respectively. The initial load levels $u_i,\forall i$ are randomly selected from $U(0\text{kW},50\text{kW})$. As for the constraint, each coefficient $\gamma_i$ is independently sampled from $U(0.03,0.15)$, and the desired load level is set to 
$D=(p_c(\mathbf{0})-1500) \text{kW}$.  
In all tests, we use the smoothing radius of $r_k=\min\{5/(k+1)^{1.1},10^{-3}\}$, and run each test 20 times with different initial points. We test all algorithms with multiple choices of step sizes, and the simulation results with the best step sizes are shown in Fig~\ref{figure convex performance} and Fig~\ref{figure convex ZOCEG}, where the relative error is computed by $(\phi_0(\mathbf{x}_k)-\phi_0^*)/\phi_0^*$. Here, $\phi_0^*$ is the optimal function value derived by an oracle solver with explicit formulation of problem \eqref{eq: constrained_problem}. The dark lines represent average values, and the shaded areas represent standard deviations.

As shown in Fig~\ref{figure convex performance}, ZOEG can converge with proper step sizes and achieve satisfactory performance in both objective optimality and constraint violation. All tested algorithms benefit from diminishing step sizes, which consistently improve the convergence behaviors. Among them, SZO-ConEX enjoys faster convergence in constraint violation but performs worse in relative error. In contrast, ZOEG and ZOGDA achieve similar and competitive performance with faster convergence in relative error.
The comparable performance of ZOEG and ZOGDA aligns with the iteration complexity proportional to $\epsilon^{-2}$ of stochastic GDA \cite{nemirovski2009robust}, which has a close relation to ZOGDA.

Figure~\ref{figure convex ZOCEG} illustrates the performance of ZOBCEG across varying BSs. Overall, ZOBCEG achieves faster convergence and smaller relative errors compared to the above algorithms. Convergence rates improve with larger BSs. Notably, ZOCEG (i.e., ZOBCEG with BS=100) converges within a few dozen iterations, aligning with our theoretical analysis and confirming improved convergence rates and computational efficiency. While ZOBCEG's theoretical complexity bounds are unfavorable, it demonstrates strong practical performance. Furthermore, Table~\ref{tab: converge_iteration} reports the average number of function evaluations required to reach specific target errors. Although ZOCEG exhibits a fast convergence rate, its requirement of $2d$ function evaluations per iteration incurs a higher single-step computational cost. In contrast, ZOBCEG with an appropriate BS offers a favorable trade-off between evaluations per iteration and total computational complexity, which can potentially reduce overall cost.



\begin{table}[htbp]
  \centering
  \caption{Average number of function evaluations required to achieve specific errors under different block sizes.}
    \begin{tabular}{ccccccc}  
    \toprule
    \multirow{2}[4]{*}{BS} & \multicolumn{3}{c}{Relative error} & \multicolumn{3}{c}{Constraint violation} \\
\cmidrule{2-7}          & 5\%   & 1\%   & 0.10\% & 5     & 1     & 0.1 \\
    \midrule
    1     & 2460.6 & 4247.1 & 5664.9 & 210.6 & 359.7 & 1309.2 \\
    5     & 905.8 & 1479.1 & 1786.4 & 183.4 & 466.2 & 1488.9 \\
    100   & 581.4 & 1458.6 & 2723.4 & 2152.2 & 2876.4 & 4324.8 \\
    \bottomrule
    \end{tabular}%
  \label{tab: converge_iteration}%
\end{table}%

\subsection{Nonconvex Case}
For the nonconvex case, we consider the demand-side energy optimization at a feeder of the 141-bus distribution network. The objective function includes the penalty for load and voltage deviation, along with the local load management costs: $\phi_0(\mathbf{x})=20\cdot(p_c(\mathbf{x})-D)^2+20\cdot\sum_{j=1}^{141}\left(\max\{v_j(\mathbf{x})-\bar{v},0\}^2+\max\{\underline{v}- v_j(\mathbf{x}),0\}^2\right)+\sum_{i}^{168}h_i(x_i)$. The total power load $p_c(\mathbf{x})$ and voltage magnitude $v_j$ at bus $j$ capture the nonlinear and nonconvex nature of AC power flows and serve as black boxes. $\mathbf{x}$ includes active and reactive power loads across all buses. and the function $h_i$ and the set $\mathcal{X}$ are generated as in the convex setting. We evaluate the performance of ZOEG under CS and DS, and ZOBCEG with various BS under CS, and adopt the technique in \cite{jin2023zeroth} to ensure the feasibility when perturbing $\mathbf{x}$.

In this case,  we set parameters $\bar{v}=1.04$ p.u., $\underline{v}=0.96$ p.u., and $D=0.15$ p.u. $=1500kW$. For ZOEG, we set $\eta=3\times10^{-6}$, $ \delta=0.005$ for CS and $\eta_k=3\times10^{-4}/\sqrt{k+1000}$, $ \delta_k=\min\{50/(k+1),0.1\}$ for DS. The smoothing radius is $r_k=\min\{0.01/(k+4000)^{1.1},10^{-5}\}$. For ZOBCEG, we consider CS with different BS, and the smoothing radius is $r_k=\min\{0.1/k^{1.2},2\times 10^{-4}\}$. 

\begin{figure}[htb]
    \centering 
    \subfloat[]{      \includegraphics[scale=0.4]{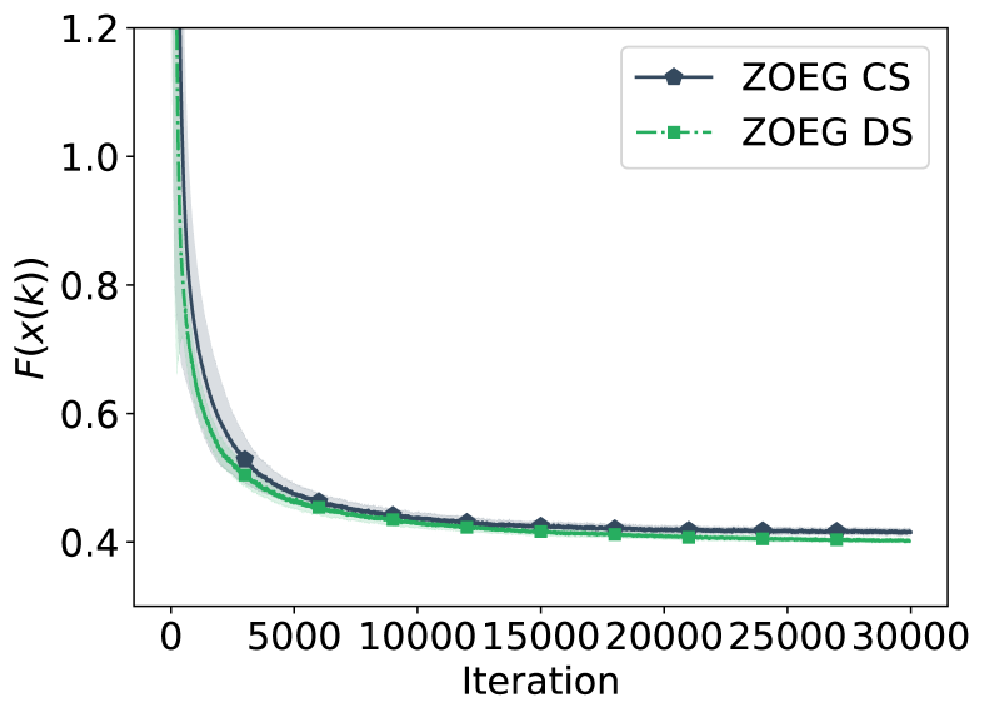}
          \label{figure nonconvex_re_ZOEG}
        }
    \hfil
    \subfloat[]{      \includegraphics[scale=0.4]{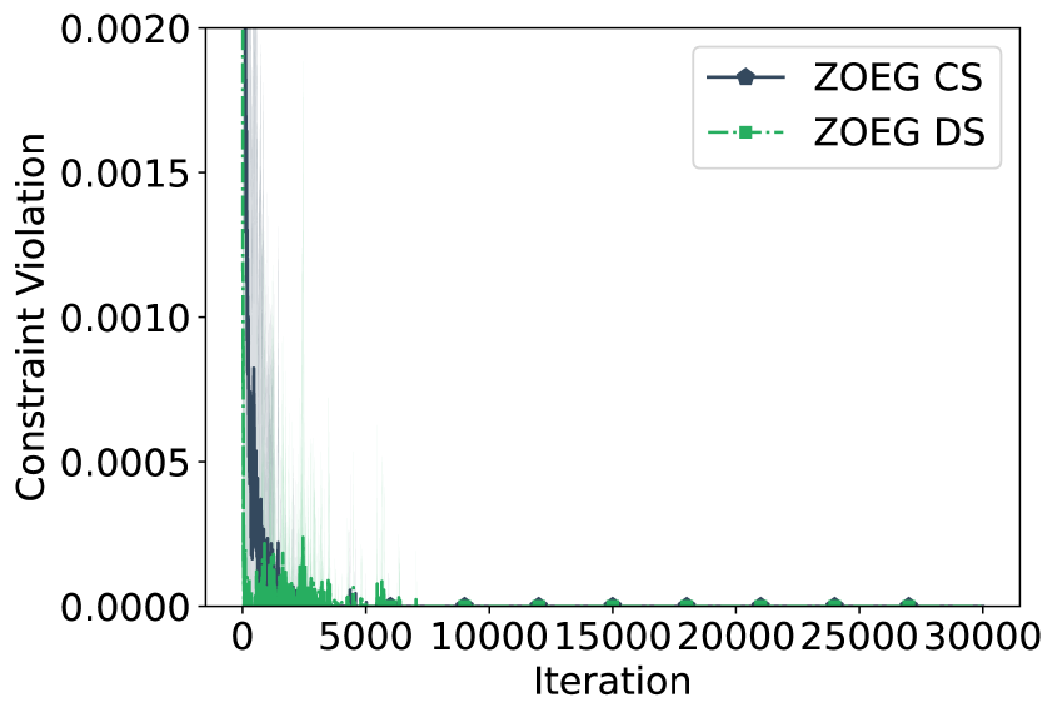}
          \label{figure nonconvex_vio_ZOEG}
    }
    \par\smallskip
    \subfloat[]{   \includegraphics[scale=0.4]{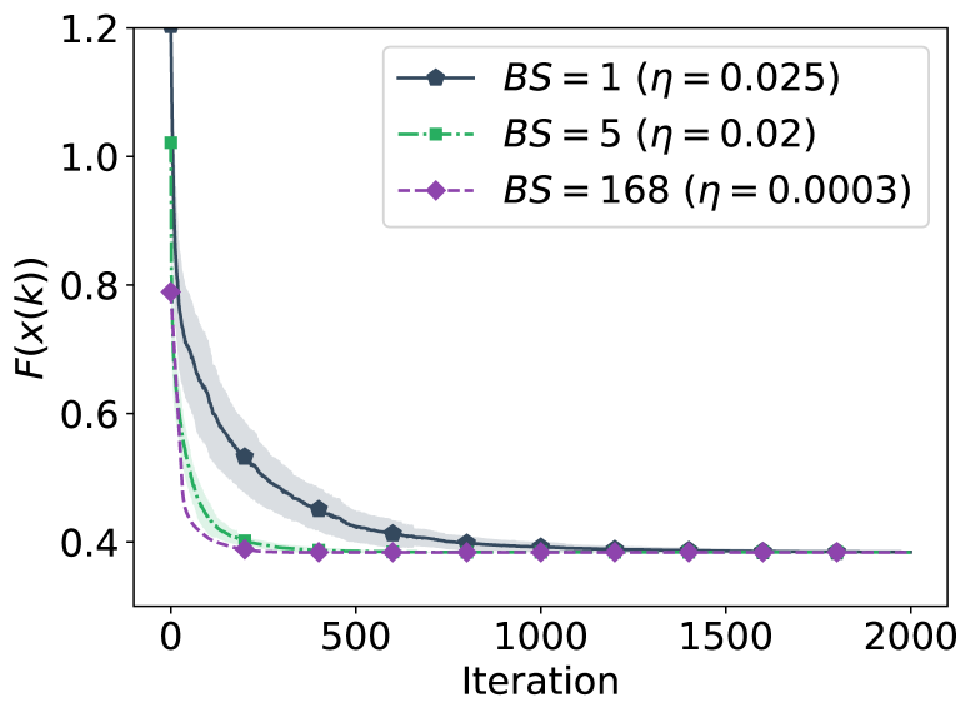}
          \label{figure nonconvex_re_ZOBCEG}
	}
    \hfil
    \subfloat[]{   \includegraphics[scale=0.4]{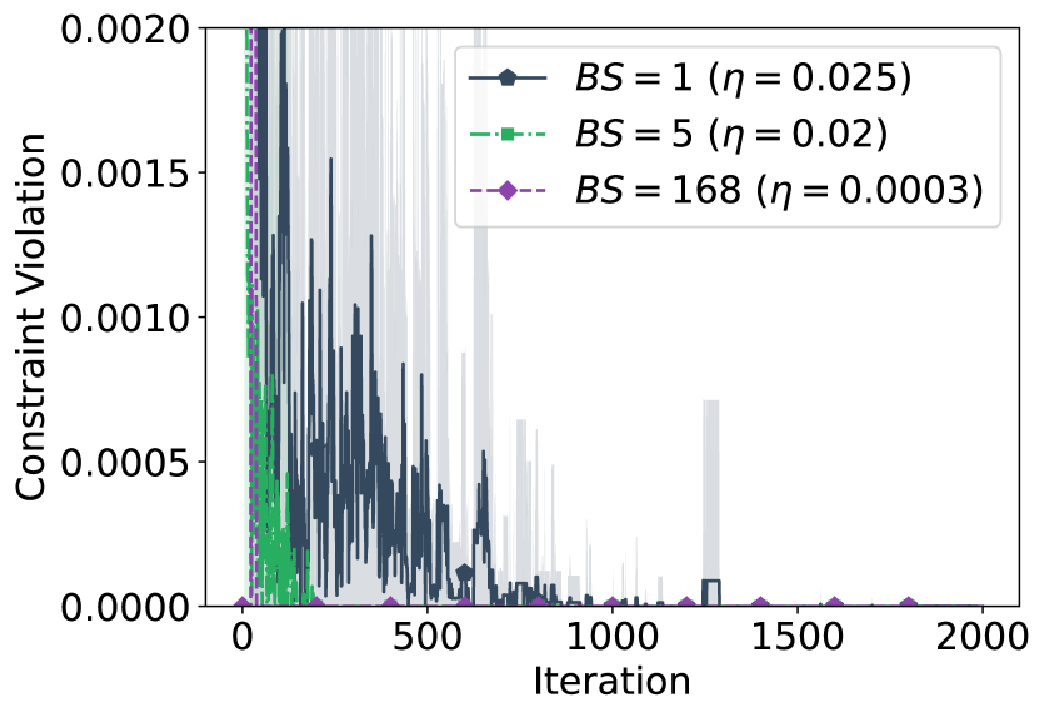}
          \label{figure nonconvex_vio_ZOBCEG}
	}
    \caption{Performance of algorithms in the nonconvex case: (a) average objective value of ZOEG, (b) average constraint violation of ZOEG, (c) average objective value of ZOBCEG, (d) average constraint violation of ZOBCEG. }
    \label{figure nonconvex performance}
    
\end{figure}

We run each case 50 times with different initial points, and the results of ZOEG and ZOBCEG are shown in Figure~\ref{figure nonconvex performance}. The proposed algorithms demonstrate excellent convergence performance even in nonconvex settings, highlighting their effectiveness in addressing more complex optimization problems. Moreover, despite having similar theoretical complexity bounds in the convex case, ZOBCEG exhibits a faster empirical convergence rate than ZOEG, and the convergence performance of ZOBCEG improves as the BS increases.

\section{Conclusions}
This work addresses a general constrained optimization problem characterized by black-box objective and constraint functions. To tackle the challenge of lacking analytical formulations, we reformulate it as a min-max problem and introduce two novel zeroth-order algorithms: the zeroth-order extra-gradient (ZOEG), leveraging a 2-point gradient estimator, and the zeroth-order coordinate extra-gradient algorithm (ZOCEG), based on a $2d$-point coordinate gradient estimator. ZOCEG is further extended to a block-coordinate variant, the zeroth-order block coordinate extra-gradient (ZOBCEG), which updates randomly selected coordinate blocks. Theoretical analysis establishes finite-sample convergence guarantees for both ZOEG and ZOCEG in convex settings regarding both optimality gap and constraint violation. Notably, ZOEG improves dependence on the problem dimension, while ZOCEG achieves the best-known convergence rate. Numerical experiments validate the efficiency of these algorithms in both convex and nonconvex scenarios.

There are several promising directions for future research: 1) Extending the algorithms to handle stochastic settings where black-box observations are corrupted by noise; 2) Improving the theoretical complexity analysis of ZOBCEG, which currently lags behind its empirical performance; 3) Establishing the theoretical guarantees for the proposed methods in nonconvex settings.


\bibliographystyle{unsrtnat}
\bibliography{Reference}

\begin{thebibliography}{42}
\providecommand{\natexlab}[1]{#1}
\providecommand{\url}[1]{\texttt{#1}}
\expandafter\ifx\csname urlstyle\endcsname\relax
  \providecommand{\doi}[1]{doi: #1}\else
  \providecommand{\doi}{doi: \begingroup \urlstyle{rm}\Url}\fi

\bibitem[Guo et~al.(2025)Guo, Jiang, Ferrari-Trecate, and Kamgarpour]{guo2025safe}
Baiwei Guo, Yuning Jiang, Giancarlo Ferrari-Trecate, and Maryam Kamgarpour.
\newblock Safe zeroth-order optimization using quadratic local approximations.
\newblock \emph{Automatica}, 174:\penalty0 112141, 2025.

\bibitem[Fu et~al.(2015)]{fu2015handbook}
Michael~C Fu et~al.
\newblock \emph{Handbook of simulation optimization}, volume 216.
\newblock Springer, 2015.

\bibitem[Nguyen and Balasubramanian(2023)]{nguyen2023stochastic}
Anthony Nguyen and Krishnakumar Balasubramanian.
\newblock Stochastic zeroth-order functional constrained optimization: Oracle complexity and applications.
\newblock \emph{INFORMS Journal on Optimization}, 5\penalty0 (3):\penalty0 256--272, 2023.

\bibitem[Chen et~al.(2021)Chen, Poveda, and Li]{chen2021modelfreeoptimalvoltagecontrol}
Xin Chen, Jorge~I. Poveda, and Na~Li.
\newblock Model-free optimal voltage control via continuous-time zeroth-order methods, 2021.
\newblock URL \url{https://arxiv.org/abs/2103.14703}.

\bibitem[Jin et~al.(2023)Jin, Tang, and Song]{jin2023zeroth}
Ruiyang Jin, Yujie Tang, and Jie Song.
\newblock Zeroth-order feedback-based optimization for distributed demand response.
\newblock \emph{arXiv preprint arXiv:2311.00372}, 2023.

\bibitem[Hooke and Jeeves(1961)]{hooke1961direct}
Robert Hooke and Terry~A Jeeves.
\newblock ``direct search''solution of numerical and statistical problems.
\newblock \emph{Journal of the ACM (JACM)}, 8\penalty0 (2):\penalty0 212--229, 1961.

\bibitem[Nelder and Mead(1965)]{nelder1965simplex}
John~A Nelder and Roger Mead.
\newblock A simplex method for function minimization.
\newblock \emph{The Computer Journal}, 7\penalty0 (4):\penalty0 308--313, 1965.

\bibitem[Fermi(1952)]{fermi1952numerical}
Enrico Fermi.
\newblock Numerical solution of a minimum problem.
\newblock Technical report, Los Alamos Scientific Lab., Los Alamos, NM, 1952.

\bibitem[Torczon(1997)]{torczon1997convergence}
Virginia Torczon.
\newblock On the convergence of pattern search algorithms.
\newblock \emph{SIAM Journal on Optimization}, 7\penalty0 (1):\penalty0 1--25, 1997.

\bibitem[{Lewis, Robert Michael and Torczon, Virginia}(1999)]{lewis1999pattern}
{Lewis, Robert Michael and Torczon, Virginia}.
\newblock Pattern search algorithms for bound constrained minimization.
\newblock \emph{SIAM Journal on Optimization}, 9\penalty0 (4):\penalty0 1082--1099, 1999.

\bibitem[Audet and Dennis~Jr(2006)]{audet2006mesh}
Charles Audet and John~E Dennis~Jr.
\newblock Mesh adaptive direct search algorithms for constrained optimization.
\newblock \emph{SIAM Journal on Optimization}, 17\penalty0 (1):\penalty0 188--217, 2006.

\bibitem[Conn et~al.(1991)Conn, Gould, and Toint]{conn1991globally}
Andrew~R Conn, Nicholas~IM Gould, and Philippe Toint.
\newblock A globally convergent augmented lagrangian algorithm for optimization with general constraints and simple bounds.
\newblock \emph{SIAM Journal on Numerical Analysis}, 28\penalty0 (2):\penalty0 545--572, 1991.

\bibitem[Ungredda and Branke(2024)]{ungredda2024bayesian}
Juan Ungredda and Juergen Branke.
\newblock Bayesian optimisation for constrained problems.
\newblock \emph{ACM Transactions on Modeling and Computer Simulation}, 34\penalty0 (2):\penalty0 1--26, 2024.

\bibitem[Regis(2011)]{regis2011stochastic}
Rommel~G Regis.
\newblock Stochastic radial basis function algorithms for large-scale optimization involving expensive black-box objective and constraint functions.
\newblock \emph{Computers \& Operations Research}, 38\penalty0 (5):\penalty0 837--853, 2011.

\bibitem[Regis(2020)]{regis2020survey}
Rommel~G Regis.
\newblock A survey of surrogate approaches for expensive constrained black-box optimization.
\newblock In \emph{Optimization of Complex Systems: Theory, Models, Algorithms and Applications}, pages 37--47. Springer, 2020.

\bibitem[Dzahini et~al.(2023)Dzahini, Kokkolaras, and Le~Digabel]{dzahini2023constrained}
Kwassi~Joseph Dzahini, Michael Kokkolaras, and S{\'e}bastien Le~Digabel.
\newblock Constrained stochastic blackbox optimization using a progressive barrier and probabilistic estimates.
\newblock \emph{Mathematical Programming}, 198\penalty0 (1):\penalty0 675--732, 2023.

\bibitem[Shi et~al.(2019)Shi, Liu, Long, Wu, and Tang]{shi2019filter}
Renhe Shi, Li~Liu, Teng Long, Yufei Wu, and Yifan Tang.
\newblock Filter-based adaptive kriging method for black-box optimization problems with expensive objective and constraints.
\newblock \emph{Computer Methods in Applied Mechanics and Engineering}, 347:\penalty0 782--805, 2019.

\bibitem[Regis and Wild(2017)]{regis2017conorbit}
Rommel~G Regis and Stefan~M Wild.
\newblock Conorbit: constrained optimization by radial basis function interpolation in trust regions.
\newblock \emph{Optimization Methods and Software}, 32\penalty0 (3):\penalty0 552--580, 2017.

\bibitem[Wang et~al.(2018)Wang, Du, Balakrishnan, and Singh]{wang2018stochastic}
Yining Wang, Simon Du, Sivaraman Balakrishnan, and Aarti Singh.
\newblock Stochastic zeroth-order optimization in high dimensions.
\newblock In \emph{International Conference on Artificial Intelligence and Statistics}, pages 1356--1365. PMLR, 2018.

\bibitem[Berahas et~al.(2022)Berahas, Cao, Choromanski, and Scheinberg]{berahas2022theoretical}
Albert~S Berahas, Liyuan Cao, Krzysztof Choromanski, and Katya Scheinberg.
\newblock A theoretical and empirical comparison of gradient approximations in derivative-free optimization.
\newblock \emph{Foundations of Computational Mathematics}, 22\penalty0 (2):\penalty0 507--560, 2022.

\bibitem[Xu et~al.(2024)Xu, Wang, Shen, and Dai]{xu2024derivative}
Zi~Xu, Ziqi Wang, Jingjing Shen, and Yuhong Dai.
\newblock Derivative-free alternating projection algorithms for general nonconvex-concave minimax problems.
\newblock \emph{SIAM Journal on Optimization}, 34\penalty0 (2):\penalty0 1879--1908, 2024.

\bibitem[Liu et~al.(2025)Liu, Wang, Xiao, and Liu]{liu2025inexact}
Shang-Lin Liu, Lei Wang, Na-Chuan Xiao, and Xin Liu.
\newblock An inexact preconditioned zeroth-order proximal method for composite optimization.
\newblock \emph{Journal of the Operations Research Society of China}, pages 1--19, 2025.

\bibitem[Zhang et~al.(2024)Zhang, Zhou, Ji, Shen, and Zavlanos]{zhang2024boosting}
Yan Zhang, Yi~Zhou, Kaiyi Ji, Yi~Shen, and Michael~M Zavlanos.
\newblock Boosting one-point derivative-free online optimization via residual feedback.
\newblock \emph{IEEE Transactions on Automatic Control}, 2024.

\bibitem[Flaxman et~al.(2005)Flaxman, Kalai, and McMahan]{flaxman2005online}
Abraham~D Flaxman, Adam~Tauman Kalai, and H~Brendan McMahan.
\newblock Online convex optimization in the bandit setting: gradient descent without a gradient.
\newblock In \emph{Proceedings of the Sixteenth Annual ACM-SIAM Symposium on Discrete Algorithms}, pages 385--394, 2005.

\bibitem[Nesterov and Spokoiny(2017)]{nesterov2017random}
Yurii Nesterov and Vladimir Spokoiny.
\newblock Random gradient-free minimization of convex functions.
\newblock \emph{Foundations of Computational Mathematics}, 17\penalty0 (2):\penalty0 527--566, 2017.

\bibitem[Tang et~al.(2023)Tang, Ren, and Li]{tang2023zeroth}
Yujie Tang, Zhaolin Ren, and Na~Li.
\newblock Zeroth-order feedback optimization for cooperative multi-agent systems.
\newblock \emph{Automatica}, 148:\penalty0 110741, 2023.

\bibitem[Duchi et~al.(2015)Duchi, Jordan, Wainwright, and Wibisono]{duchi2015optimal}
John~C Duchi, Michael~I Jordan, Martin~J Wainwright, and Andre Wibisono.
\newblock Optimal rates for zero-order convex optimization: The power of two function evaluations.
\newblock \emph{IEEE Transactions on Information Theory}, 61\penalty0 (5):\penalty0 2788--2806, 2015.

\bibitem[Liu et~al.(2020{\natexlab{a}})Liu, Chen, Kailkhura, Zhang, Hero~III, and Varshney]{liu2020primer}
Sijia Liu, Pin-Yu Chen, Bhavya Kailkhura, Gaoyuan Zhang, Alfred~O Hero~III, and Pramod~K Varshney.
\newblock A primer on zeroth-order optimization in signal processing and machine learning: Principals, recent advances, and applications.
\newblock \emph{IEEE Signal Processing Magazine}, 37\penalty0 (5):\penalty0 43--54, 2020{\natexlab{a}}.

\bibitem[Ghadimi et~al.(2016)Ghadimi, Lan, and Zhang]{ghadimi2016mini}
Saeed Ghadimi, Guanghui Lan, and Hongchao Zhang.
\newblock Mini-batch stochastic approximation methods for nonconvex stochastic composite optimization.
\newblock \emph{Mathematical Programming}, 155\penalty0 (1):\penalty0 267--305, 2016.

\bibitem[Yu et~al.(2021)Yu, Ho, and Yuan]{yu2021distributed}
Zhan Yu, Daniel~WC Ho, and Deming Yuan.
\newblock Distributed randomized gradient-free mirror descent algorithm for constrained optimization.
\newblock \emph{IEEE Transactions on Automatic Control}, 67\penalty0 (2):\penalty0 957--964, 2021.

\bibitem[Chen et~al.(2020)Chen, Zhou, Yi, and Gu]{chen2020frank}
Jinghui Chen, Dongruo Zhou, Jinfeng Yi, and Quanquan Gu.
\newblock A frank-wolfe framework for efficient and effective adversarial attacks.
\newblock In \emph{Proceedings of the AAAI Conference on Artificial Intelligence}, volume~34, pages 3486--3494, 2020.

\bibitem[Ye et~al.(2025)Ye, Huang, Di, and Chang]{ye2025enhanced}
Haishan Ye, Yinghui Huang, Hao Di, and Xiangyu Chang.
\newblock An enhanced zeroth-order stochastic frank-wolfe framework for constrained finite-sum optimization.
\newblock \emph{arXiv preprint arXiv:2501.07201}, 2025.

\bibitem[Hu et~al.(2023)Hu, Zhang, and Wu]{hu2023gradient}
Chuanhao Hu, Xuan Zhang, and Qiuwei Wu.
\newblock Gradient-free accelerated event-triggered scheme for constrained network optimization in smart grids.
\newblock \emph{IEEE Transactions on Smart Grid}, 15\penalty0 (3):\penalty0 2843--2855, 2023.

\bibitem[Liu et~al.(2020{\natexlab{b}})Liu, Lu, Chen, Feng, Xu, Al-Dujaili, Hong, and O’Reilly]{liu2020min}
Sijia Liu, Songtao Lu, Xiangyi Chen, Yao Feng, Kaidi Xu, Abdullah Al-Dujaili, Mingyi Hong, and Una-May O’Reilly.
\newblock Min-max optimization without gradients: Convergence and applications to black-box evasion and poisoning attacks.
\newblock In \emph{International Conference on Machine Learning}, pages 6282--6293. PMLR, 2020{\natexlab{b}}.

\bibitem[Maheshwari et~al.(2022)Maheshwari, Chiu, Mazumdar, Sastry, and Ratliff]{maheshwari2022zeroth}
Chinmay Maheshwari, Chih-Yuan Chiu, Eric Mazumdar, Shankar Sastry, and Lillian Ratliff.
\newblock Zeroth-order methods for convex-concave min-max problems: Applications to decision-dependent risk minimization.
\newblock In \emph{International Conference on Artificial Intelligence and Statistics}, pages 6702--6734. PMLR, 2022.

\bibitem[Nemirovski(2004)]{nemirovski2004prox}
Arkadi Nemirovski.
\newblock Prox-method with rate of convergence o (1/t) for variational inequalities with lipschitz continuous monotone operators and smooth convex-concave saddle point problems.
\newblock \emph{SIAM Journal on Optimization}, 15\penalty0 (1):\penalty0 229--251, 2004.

\bibitem[Nedi{\'c} and Ozdaglar(2009)]{nedic2009subgradient}
Angelia Nedi{\'c} and Asuman Ozdaglar.
\newblock Subgradient methods for saddle-point problems.
\newblock \emph{Journal of Optimization Theory and Applications}, 142:\penalty0 205--228, 2009.

\bibitem[Malik et~al.(2020)Malik, Pananjady, Bhatia, Khamaru, Bartlett, and Wainwright]{malik2020derivative}
Dhruv Malik, Ashwin Pananjady, Kush Bhatia, Koulik Khamaru, Peter~L Bartlett, and Martin~J Wainwright.
\newblock Derivative-free methods for policy optimization: Guarantees for linear quadratic systems.
\newblock \emph{Journal of Machine Learning Research}, 21\penalty0 (21):\penalty0 1--51, 2020.

\bibitem[Mokhtari et~al.(2020)Mokhtari, Ozdaglar, and Pattathil]{mokhtari2020unified}
Aryan Mokhtari, Asuman Ozdaglar, and Sarath Pattathil.
\newblock A unified analysis of extra-gradient and optimistic gradient methods for saddle point problems: Proximal point approach.
\newblock In \emph{International Conference on Artificial Intelligence and Statistics}, pages 1497--1507. PMLR, 2020.

\bibitem[He et~al.(2022)He, Lu, Guan, Kang, and Shi]{he2022zeroth}
Pengcheng He, Siyuan Lu, Xin Guan, Yibin Kang, and Qingjiang Shi.
\newblock A zeroth-order block coordinate gradient descent method for cellular network optimization.
\newblock In \emph{2022 International Symposium on Wireless Communication Systems (ISWCS)}, pages 1--6. IEEE, 2022.

\bibitem[Boob et~al.(2023)Boob, Deng, and Lan]{boob2023stochastic}
Digvijay Boob, Qi~Deng, and Guanghui Lan.
\newblock Stochastic first-order methods for convex and nonconvex functional constrained optimization.
\newblock \emph{Mathematical Programming}, 197\penalty0 (1):\penalty0 215--279, 2023.

\bibitem[Nemirovski et~al.(2009)Nemirovski, Juditsky, Lan, and Shapiro]{nemirovski2009robust}
Arkadi Nemirovski, Anatoli Juditsky, Guanghui Lan, and Alexander Shapiro.
\newblock Robust stochastic approximation approach to stochastic programming.
\newblock \emph{SIAM Journal on Optimization}, 19\penalty0 (4):\penalty0 1574--1609, 2009.

\end{thebibliography}

\newpage

\begin{center}
    {\LARGE\bfseries Appendices}
\end{center}

\appendix

\section{Proof of Lemma~\ref{lem: coor_grad_estimator}}\label{apx: pf_coor_grad_est}
    By definition, we have
    \begin{align*}  
    &\ \left |G_{h,i}(\mathbf{z};r)-\nabla_{i}h(\mathbf{z})\right|\\  =&\ \left|\frac{h(\mathbf{z}+r\mathbf{e}_i)-h(\mathbf{z})}{r}-\nabla_{i} h(\mathbf{z})\right|\\
   =&\ \left|\int_0^1 \nabla_{i} h(\mathbf{z}+s\cdot r\mathbf{e}_i)-\nabla_{i} h(\mathbf{z}_k)ds \right|\\ 
  \leq&\ \int_0^1\left|\nabla_{i} h(\mathbf{z}+s\cdot r\mathbf{e}_i)-\nabla_{i} h(\mathbf{z})\right|ds\\
        \leq&\ \int_0^1 L\cdot\left|s\cdot r\mathbf{e}_i\right|ds\leq \frac{1}{2}Lr,
    \end{align*}
where the second inequality follows from the smoothness.
\subsection{Proof of Lemma~\ref{lem: bound_EG_term}}\label{pf: bound_EG_term}
Firstly, by introducing $\left\langle\eta F^{r_k}(\mathbf{z}_k^+),\mathbf{z}_k^+-\mathbf{z}\right\rangle$, we have:
\begin{align*}
    &\ \left\langle\eta F(\mathbf{z}_k^+),\mathbf{z}_k^+-\mathbf{z}\right\rangle\\
            =&\ \mathbb{E}_{\mathbf{w}_k}\left[\left\langle\eta\mathbf{g}_k^+,\mathbf{z}_k^+-\mathbf{z}\right\rangle\right]+\left\langle\eta F(\mathbf{z}_k^+)-\eta F^{r_k}(\mathbf{z}_k^+),\mathbf{z}_k^+-\mathbf{z}\right\rangle\\
           \leq &\ \mathbb{E}_{\mathbf{w}_k}\left[\left\langle\eta\mathbf{g}_k^+,\mathbf{z}_k^+-\mathbf{z}\right\rangle\right] + \eta \Vert F(\mathbf{z}_k^+)-F^{r_k}(\mathbf{z}_k^+)\Vert \Vert \mathbf{z}_k^+-\mathbf{z}\Vert\\
           \leq &\ \mathbb{E}_{\mathbf{w}_k}\left[\left\langle\eta\mathbf{g}_k^+,\mathbf{z}_k^+-\mathbf{z}\right\rangle\right] +\eta Lr_k\cdot\tilde{D},
\end{align*}
   where the first equality follows from \eqref{eq: exp_2p_grad} and the last inequality follows from \eqref{eq: diff smooth}.
   Taking expectation over $\mathbf{v}_k$ gives:
      \begin{equation}
    \begin{aligned}\label{eq: pf bound cond first order}
       &\ \mathbb{E}\left[\left\langle\eta F(\mathbf{z}_k^+),\mathbf{z}_k^+-\mathbf{z}\right\rangle\mid \mathcal{F}_k\right]\\
        \leq &\ \mathbb{E}\left[\left\langle\eta\mathbf{g}_k^+,\mathbf{z}_k^+-\mathbf{z}\right\rangle\mid \mathcal{F}_k\right]+\eta L\tilde{D}r_k.
   \end{aligned}    
   \end{equation}
   Then, we focus on the bound for the first term:
   \begin{align}\label{eq: pf bound first term}
       &\ \left\langle\eta\mathbf{g}_k^+,\mathbf{z}_k^+-\mathbf{z}\right\rangle\nonumber \\
        =&\     \left\langle\eta\mathbf{g}_k^+,\mathbf{z}_k^+-\mathbf{z}_{k+1}\right\rangle+   \left\langle\eta\mathbf{g}_k^+,\mathbf{z}_{k+1}-\mathbf{z}\right\rangle\nonumber\\
               \leq&\    \left\langle\eta\mathbf{g}_k^+,\mathbf{z}_k^+-\mathbf{z}_{k+1}\right\rangle-   \left\langle\mathbf{z}_{k+1}-\mathbf{z}_{k},\mathbf{z}_{k+1}-\mathbf{z}\right\rangle\nonumber\\
               =&\ \left\langle\eta\mathbf{g}_k^+,\mathbf{z}_k^+-\mathbf{z}_{k+1}\right\rangle-\frac{1}{2}\Vert \mathbf{z}_{k+1}-\mathbf{z}_{k} \Vert^2+H_{k}(\mathbf{z})-H_{k+1}(\mathbf{z}),
   \end{align}
    where the inequality follow from the projection property:
     \begin{equation*} 
\left\langle\eta\mathbf{g}_k^++\mathbf{z}_{k+1}-\mathbf{z}_k,\mathbf{z}_{k+1}-\mathbf{z}\right\rangle\leq 0.\label{eq:proj_k_grad}
    \end{equation*}  
    Similarly, by introducing the auxiliary term $\langle\eta\mathbf{g}_k,\mathbf{z}_k^+-\mathbf{z}_{k+1}\rangle$, we can obtain the upper bound for $\langle\eta\mathbf{g}_k^+,\mathbf{z}_k^+-\mathbf{z}_{k+1}\rangle-\frac{1}{2}\Vert \mathbf{z}_{k+1}-\mathbf{z}_{k} \Vert^2$ as follows:
    \begin{align}\label{eq: pf bound first two term}
        &\ \left\langle\eta\mathbf{g}_k^+,\mathbf{z}_k^+-\mathbf{z}_{k+1}\right\rangle-\frac{1}{2}\Vert \mathbf{z}_{k+1}-\mathbf{z}_{k} \Vert^2\nonumber\\
        =&\ \eta\left\langle\mathbf{g}_k^+-\mathbf{g}_k,\mathbf{z}_k^+-\mathbf{z}_{k+1}\right\rangle+\left\langle\eta\mathbf{g}_k,\mathbf{z}_k^+-\mathbf{z}_{k+1}\right\rangle-\frac{1}{2}\Vert \mathbf{z}_{k+1}-\mathbf{z}_{k} \Vert^2\nonumber\\
          \leq&\ \eta\left\langle\mathbf{g}_k^+-\mathbf{g}_k,\mathbf{z}_k^+-\mathbf{z}_{k+1}\right\rangle - \left\langle\mathbf{z}_k^+-\mathbf{z}_k,\mathbf{z}_k^+-\mathbf{z}_{k+1}\right\rangle-\frac{1}{2}\Vert \mathbf{z}_{k+1}-\mathbf{z}_{k} \Vert^2\nonumber\\     
            =&\ \eta\left\langle\mathbf{g}_k^+-\mathbf{g}_k,\mathbf{z}_k^+-\mathbf{z}_{k+1}\right\rangle -\frac{1}{2}\Vert \mathbf{z}_{k}^+-\mathbf{z}_{k} \Vert^2-\frac{1}{2}\Vert \mathbf{z}_{k}^+-\mathbf{z}_{k+1} \Vert^2\nonumber\\
            \leq &\ \eta\Vert\mathbf{g}_k-\mathbf{g}_{k}^+\Vert\Vert\mathbf{z}_k^+-\mathbf{z}_{k+1}\Vert\nonumber\\
            \leq&\ \eta^2 \Vert\mathbf{g}_k-\mathbf{g}_{k}^+\Vert^2,
    \end{align}
    where the first inequality uses the projection property: 
    \begin{equation*}
    \left\langle\eta\mathbf{g}_k+\mathbf{z}_k^+-\mathbf{z}_k,\mathbf{z}_k^+-\mathbf{z}\right\rangle\leq 0,\label{eq:proj_plus_grad}
    \end{equation*}  
    and the final bound follows from Cauchy-Schwarz inequality and the contraction property of the projection operator.
    
    Since $$\mathbf{g}_k-\mathbf{g}_{k}^+=\mathbf{g}_k-F(\mathbf{z}_k)+F(\mathbf{z}_k)-F(\mathbf{z}_k^+)+F(\mathbf{z}_k^+)-\mathbf{g}_{k}^+,$$
    we can bound the conditional expectation of $\Vert\mathbf{g}_k-\mathbf{g}_{k}^+\Vert^2$:
    \begin{align*}
&\ \mathbb{E}\left[\Vert\mathbf{g}_k-\mathbf{g}_{k}^+\Vert^2\mid\mathcal{F}_k\right]\\
        \leq &\ 3\cdot\mathbb{E}\left[\Vert\mathbf{g}_k-F(\mathbf{z}_k)\Vert^2+\Vert\mathbf{g}_k^+-F(\mathbf{z}_k^+)\Vert^2\mid\mathcal{F}_k\right]+3\cdot\mathbb{E}\left[\Vert F(\mathbf{z}_k)-F(\mathbf{z}_k^+)\Vert^2\mid\mathcal{F}_k\right]\\
        \leq&\ 3\cdot\mathbb{E}\left[\Vert\mathbf{g}_k-F(\mathbf{z}_k)\Vert^2\mid\mathcal{F}_k\right]+3\cdot\mathbb{E}\left[\Vert\mathbf{g}_k^+-F(\mathbf{z}_k^+)\Vert^2\mid\mathcal{F}_k\right]+12G^2.
    \end{align*}
    Note that combining the bound for variance of gradient estimator \eqref{eq: diff smooth} with \eqref{eq: var_2p_grad} yields the following bound:
    \begin{equation*}
        \begin{aligned}
            &\mathbb{E}\left[\Vert\mathbf{g}_k-F(\mathbf{z}_k)\Vert^2\mid\mathcal{F}_k\right]\\
            \leq&2\cdot\mathbb{E}\left[\Vert\mathbf{g}_k-F^{r_k}(\mathbf{z}_k)\Vert^2+\Vert F^{r_k}(\mathbf{z}_k)-F(\mathbf{z}_k)\Vert^2\mid\mathcal{F}_k\right]\\
            \leq &\frac{18d^2G^2}{d+2}+\frac{3d^2 L^2 r_k^2}{2}+2L^2r_k^2,
        \end{aligned}
    \end{equation*}
    and the bound on $\mathbb{E}\left[\Vert\mathbf{g}_k^+-F(\mathbf{z}_k^+)\Vert^2\mid\mathcal{F}_k\right]$ can be derived similarly. Then, $\mathbb{E}\left[\Vert\mathbf{g}_k-\mathbf{g}_{k}^+\Vert^2\mid\mathcal{F}_k\right]$ is bounded as follows:
    \begin{equation}\label{eq: pf bound cond var g}
       \mathbb{E}\left[\Vert\mathbf{g}_k-\mathbf{g}_{k}^+\Vert^2\mid\mathcal{F}_k\right]\leq \hat{\Delta}_k,
    \end{equation}
    where $\hat{\Delta}_k=\frac{108d^2G^2}{d+2}+9d^2 L^2 r_k^2+12L^2r_k^2+12G^2$. The proof is completed by combining the bounds established in \eqref{eq: pf bound cond var g}, \eqref{eq: pf bound first two term}, \eqref{eq: pf bound first term}, and \eqref{eq: pf bound cond first order}.

\section{Proof of Lemma~\ref{lem: bound ZOCEG term}}\label{pf: bound ZOCEG term}
For notation simplicity, we denote $H_k^{x}(\mathbf{x})=\frac{1}{2}\Vert\mathbf{x}_k-\mathbf{x}\Vert^2$, and $H_k^{y}(\mathbf{y})=\frac{1}{2}\Vert\mathbf{y}_k-\mathbf{y}\Vert^2$.
    First,  we deal with $\left\langle\eta\nabla_{\mathbf{x}}f(\mathbf{z}_k^+),\mathbf{x}_k^{+}-\mathbf{x}\right\rangle$ by introducing $\eta \left\langle \hat{\mathbf{g}}_k^{x,+},\mathbf{x}_k^{+}-\mathbf{x}\right\rangle$: 
\begin{align} \label{eq: full_grad_bound_ZOCEG_term}  
&\  \left\langle\eta\nabla_{\mathbf{x}}f(\mathbf{z}_k^+),\mathbf{x}_k^{+}-\mathbf{x}\right\rangle\nonumber\\
    =&\ \eta \left\langle \hat{\mathbf{g}}_k^{x,+}+\nabla_{\mathbf{x}}f(\mathbf{z}_k^+)-\hat{\mathbf{g}}_k^{x,+},\mathbf{x}_k^{+}-\mathbf{x}\right\rangle\nonumber\\
    \leq&\ \eta \left\langle \hat{\mathbf{g}}_k^{x,+},\mathbf{x}_k^{+}-\mathbf{x}\right\rangle+\frac{1}{2}\eta Lr_k\sqrt{d_x}\tilde{D},
\end{align}
where the last inequality follows from the bounded bias of coordinate gradient estimator: $\Vert \nabla_{\mathbf{x}}f(\mathbf{z}_k^+)-\hat{\mathbf{g}}_k^{x,+}\Vert\leq\frac{\sqrt{d_x}}{2}Lr_k$, which directly follows from \eqref{eq: coor_grad_est_bias}.
    
    We then focus on the bound on the first term $\eta \left\langle \hat{\mathbf{g}}_k^{x,+},\mathbf{x}_k^{+}-\mathbf{x}\right\rangle$ as follows:
        \begin{align}\label{eq: full_grad_bound_first_term}
          &\ \eta \left\langle \hat{\mathbf{g}}_k^{x,+},\mathbf{x}_k^{+}-\mathbf{x}\right\rangle\nonumber\\ 
          =&\ \eta \left\langle \hat{\mathbf{g}}_k^{x,+},\mathbf{x}_k^{+}-\mathbf{x}_{k+1}\right\rangle+\eta \left\langle \hat{\mathbf{g}}_k^{x,+},\mathbf{x}_{k+1}-\mathbf{x}\right\rangle\nonumber\\
        \leq &\ \eta \left\langle \hat{\mathbf{g}}_k^{x,+},\mathbf{x}_k^{+}-\mathbf{x}_{k+1}\right\rangle-\left\langle \mathbf{x}_{k+1}-\mathbf{x}_k,\mathbf{x}_{k+1}-\mathbf{x}\right\rangle\nonumber\\
          =&\ \eta \left\langle \hat{\mathbf{g}}_k^{x,+}-\hat{\mathbf{g}}_k^{x},\mathbf{x}_k^{+}-\mathbf{x}_{k+1}\right\rangle+\eta \left\langle \hat{\mathbf{g}}_k^{x},\mathbf{x}_k^{+}-\mathbf{x}_{k+1}\right\rangle-\left\langle \mathbf{x}_{k+1}-\mathbf{x}_k,\mathbf{x}_{k+1}-\mathbf{x}\right\rangle\nonumber\\
           \leq&\ \eta \left\langle \hat{\mathbf{g}}_k^{x,+}-\hat{\mathbf{g}}_k^{x},\mathbf{x}_k^{+}-\mathbf{x}_{k+1}\right\rangle
           -\left\langle \mathbf{x}_{k}^+-\mathbf{x}_k,\mathbf{x}_{k}^+-\mathbf{x}_{k+1}\right\rangle -\left\langle \mathbf{x}_{k+1}-\mathbf{x}_k,\mathbf{x}_{k+1}-\mathbf{x}\right\rangle,
        \end{align}
    where the two inequalities follow from the properties of projection:
    \begin{equation*}
         \left\langle \eta\hat{\mathbf{g}}_k^{x,+}+\mathbf{x}_{k+1}-\mathbf{x}_{k},\mathbf{x}_{k+1}-\mathbf{x}\right\rangle\leq 0,
    \end{equation*}
     \begin{equation*}
         \left\langle \eta\hat{\mathbf{g}}_k^{x}+\mathbf{x}_{k}^+-\mathbf{x}_{k},\mathbf{x}_{k}^+-\mathbf{x}\right\rangle\leq 0.
    \end{equation*}
By performing a simple algebraic transformation, we have:
\begin{equation}
\begin{aligned}\label{eq: full_grad_alge}
&\ \left\langle \mathbf{x}_{k}-\mathbf{x}_k^+,\mathbf{x}_{k}^+-\mathbf{x}_{k+1}\right\rangle+\left\langle \mathbf{x}_{k}-\mathbf{x}_{k+1},\mathbf{x}_{k+1}-\mathbf{x}\right\rangle\\
=&\ H_{k}^{x}(\mathbf{x})-H_{k+1}^{x}(\mathbf{x})-\frac{1}{2}\Vert\mathbf{x}_{k+1}-\mathbf{x}_k^+\Vert^2-\frac{1}{2}\Vert\mathbf{x}_{k}^+-\mathbf{x}_k\Vert^2.
\end{aligned}
\end{equation}
Besides,  $\left\langle \hat{\mathbf{g}}_k^{x,+}-\hat{\mathbf{g}}_k^{x},\mathbf{x}_k^{+}-\mathbf{x}_{k+1}\right\rangle$ can be bounded as follows:
\begin{align}\label{eq: full_grad_FD}
        &\ \left\langle \hat{\mathbf{g}}_k^{x,+}-\hat{\mathbf{g}}_k^{x},\mathbf{x}_k^{+}-\mathbf{x}_{k+1}\right\rangle\nonumber\\
        =&\ \left\langle \hat{\mathbf{g}}_k^{x,+}-\nabla_{\mathbf{x}}f(\mathbf{z}_k^+)+\nabla_{\mathbf{x}}f(\mathbf{z}_k)-
\hat{\mathbf{g}}_k^{x},\mathbf{x}_k^{+}-\mathbf{x}_{k+1}\right\rangle +\left\langle \nabla_{\mathbf{x}}f(\mathbf{z}_k^+)-\nabla_{\mathbf{x}}f(\mathbf{z}_k),\mathbf{x}_k^{+}-\mathbf{x}_{k+1}\right\rangle\nonumber\\
        \leq&\  Lr_k\sqrt{d_x}\tilde{D}+L\Vert\mathbf{x}_k^+-\mathbf{x}_k\Vert\Vert\mathbf{x}_k^+-\mathbf{x}_{k+1}\Vert+L\Vert\mathbf{y}_k^+-\mathbf{y}_k\Vert\Vert\mathbf{x}_k^+-\mathbf{x}_{k+1}\Vert,
\end{align}
where the last step follows from Lemma \ref{lem: coor_grad_estimator} and the smoothness of $f$.
By the condition $\eta L\leq\frac{1}{2}$, substituting \eqref{eq: full_grad_alge}, \eqref{eq: full_grad_FD}, and \eqref{eq: full_grad_bound_first_term} into \eqref{eq: full_grad_bound_ZOCEG_term} gives:
    \begin{align}\label{eq: full_grad_ZOCEG_x_term}
        &\ \left\langle\eta\nabla_{\mathbf{x}}f(\mathbf{z}_k^+),\mathbf{x}_k^+-\mathbf{x}\right\rangle\nonumber\\ 
        \leq&\ \frac{3}{2}\eta Lr_k\sqrt{d_x}\tilde{D}+\eta L\Vert\mathbf{x}_k^+-\mathbf{x}_k\Vert\Vert\mathbf{x}_k^+-\mathbf{x}_{k+1}\Vert+H_{k}^{x}(\mathbf{x})-H_{k+1}^{x}(\mathbf{x})-\frac{1}{2}\Vert\mathbf{x}_{k+1}-\mathbf{x}_k^+\Vert^2\nonumber\\
        &\ -\frac{1}{2}\Vert\mathbf{x}_{k}^+-\mathbf{x}_k\Vert^2+\eta L\Vert\mathbf{y}_k^+-\mathbf{y}_k\Vert\Vert\mathbf{x}_k^+-\mathbf{x}_{k+1}\Vert\nonumber\\
        \leq&\ H_{k}^{x}(\mathbf{x})-H_{k+1}^{x}(\mathbf{x})+\frac{3}{2}\eta Lr_k\sqrt{d_x}\tilde{D}-\frac{1}{4}\Vert\mathbf{x}_{k}^+-\mathbf{x}_k\Vert^2+\frac{1}{4}\Vert\mathbf{y}_k^+-\mathbf{y}_k\Vert^2.
\end{align}
where in the last step we used $a\cdot b \leq \frac{1}{2}(a^2+b^2), \forall a,b\in\mathbb{R}$. Similarly, we have:
\begin{equation}
    \label{eq: full_grad_ZOCEG_y_term}
         -\left\langle\eta\nabla_{\mathbf{y}}f(\mathbf{z}_k^+),\mathbf{y}_k^+-\mathbf{y}\right\rangle 
         \leq H_{k}^{y}(\mathbf{y})-H_{k+1}^{y}(\mathbf{y})+\frac{3}{2}\eta Lr_k\sqrt{d_y}\tilde{D}-\frac{1}{4}\Vert\mathbf{y}_{k}^+-\mathbf{y}_k\Vert^2+\frac{1}{4}\Vert\mathbf{x}_k^+-\mathbf{x}_k\Vert^2.
\end{equation}
    Adding \eqref{eq: full_grad_ZOCEG_x_term} and  \eqref{eq: full_grad_ZOCEG_y_term} gives:
    \begin{equation*}
        \left\langle \eta F(\mathbf{z}_k^+),\mathbf{z}_k^+-\mathbf{z}\right\rangle
        \leq H_k(\mathbf{z})-H_{k+1}(\mathbf{z})+\frac{3}{2}\eta Lr_k(\sqrt{d_x}+\sqrt{d_y})\tilde{D}.
    \end{equation*}

\section{Proof of Theorem~\ref{thm: ZOBCEG convergence}}\label{apx: pf_ZOBCEG}
Similarly, We start by defining a filtration $\mathcal{F}_k^b:=\sigma\left(\mathbf{z}_0, \mathcal{I}_0^x,\mathcal{I}_0^y,\mathcal{J}_0^x,\mathcal{J}_0^y,\mathbf{z}_1, ...,\mathbf{z}_k\right)$. Given $\mathbf{z}_k$, $\mathbf{z}_k^+$,
$\langle\nabla _{\mathbf{x}}f(\mathbf{z}_k^+),\allowbreak\mathbf{x}_k^+-\mathbf{x}\rangle$ is bounded as follows:
\begin{align*}
    &\ \left\langle\nabla _{\mathbf{x}}f(\mathbf{z}_k^+),\mathbf{x}_k^+-\mathbf{x}\right\rangle\nonumber\\
    \leq&\ \left\langle\mathbb{E}_{\mathcal{J}_k^x}\left[\frac{d_x}{\tau_x}\tilde{\mathbf{g}}_k^{x,+}\right],\mathbf{x}_k^+-\mathbf{x}\right\rangle +\left\Vert\nabla _{\mathbf{x}}f(\mathbf{z}_k^+)-\mathbb{E}_{\mathcal{J}_k^x}\left[\frac{d_x}{\tau_x}\tilde{\mathbf{g}}_k^{x,+}\right]\right\Vert\cdot\Vert\mathbf{x}_k^+-\mathbf{x}\Vert\nonumber\\
    \leq&\ \left\langle\mathbb{E}_{\mathcal{J}_k^x}\left[\frac{d_x}{\tau_x}\tilde{\mathbf{g}}_k^{x,+}\right],\mathbf{x}_k^+-\mathbf{x}\right\rangle+\frac{Lr_k\sqrt{d_x}\tilde{D}}{2} . 
\end{align*}
Then, taking expectation over $\mathcal{I}_k^x,\mathcal{I}_k^y$ gives:
\begin{equation}
   \begin{aligned}
    \label{eq: block_grad_bound_ZOCEG_term}&\ \mathbb{E}\left[\left.\left\langle\nabla _{\mathbf{x}}f(\mathbf{z}_k^+),\mathbf{x}_k^+-\mathbf{x}\right\rangle\right| \mathcal{F}_k^b\right] \leq \mathbb{E}\left[\left.\left\langle\frac{d_x}{\tau_x}\tilde{\mathbf{g}}_k^{x,+},\mathbf{x}_k^+-\mathbf{x}\right\rangle\right| \mathcal{F}_k^b\right]+\frac{Lr_k\sqrt{d_x}\tilde{D}}{2} . 
\end{aligned} 
\end{equation}
Specifically, for $\left\langle\tilde{\mathbf{g}}_k^{x,+},\mathbf{x}_k^+-\mathbf{x}\right\rangle$, we have:
\begin{align}
    \label{eq: block_grad_bound_first_term}&\ \left\langle\tilde{\mathbf{g}}_k^{x,+},\mathbf{x}_k^+-\mathbf{x}\right\rangle\nonumber\\
    \leq&\ \left\langle\tilde{\mathbf{g}}_k^{x,+},\mathbf{x}_k^+-\mathbf{x}_{k+1}\right\rangle-\frac{1}{\eta}\left\langle\mathbf{x}_{k+1}-\mathbf{x}_k,\mathbf{x}_{k+1}-\mathbf{x}\right\rangle\nonumber\\
    \leq&\ \left\langle\tilde{\mathbf{g}}_k^{x,+},\mathbf{x}_k^+-\mathbf{x}_{k+1}\right\rangle+\frac{H_k^x(\mathbf{x})-H_{k+1}^x(\mathbf{x})}{\eta}-\frac{1}{2\eta}\Vert\mathbf{x}_{k}-\mathbf{x}_{k+1}\Vert^2\nonumber\\
    \leq&\ \left\langle\tilde{\mathbf{g}}_k^{x,+}-\tilde{\mathbf{g}}_k^{x},\mathbf{x}_k^+-\mathbf{x}_{k+1}\right\rangle-\frac{1}{\eta}\left\langle\mathbf{x}_k^+-\mathbf{x}_k,\mathbf{x}_k^+-\mathbf{x}_{k+1}\right\rangle+\frac{H_k^x(\mathbf{x})-H_{k+1}^x(\mathbf{x})}{\eta}-\frac{1}{2\eta}\Vert\mathbf{x}_{k}-\mathbf{x}_{k+1}\Vert^2\nonumber\\
    \leq&\ \left\langle\tilde{\mathbf{g}}_k^{x,+}-\tilde{\mathbf{g}}_k^{x},\mathbf{x}_k^+-\mathbf{x}_{k+1}\right\rangle+\frac{H_k^x(\mathbf{x})-H_{k+1}^x(\mathbf{x})}{\eta}\nonumber\\
    \leq&\ \eta \Vert\tilde{\mathbf{g}}_k^{x,+}-\tilde{\mathbf{g}}_k^{x}\Vert^2+\frac{H_k^x(\mathbf{x})-H_{k+1}^x(\mathbf{x})}{\eta},   
\end{align}
where the first, third, and last inequalities follow from the projection property, while the second and fourth arise from straightforward algebraic transformations
.
The conditional expectation of $\Vert\tilde{\mathbf{g}}_k^{x,+}-\tilde{\mathbf{g}}_k^{x}\Vert^2$ can be bounded as follows:
\begin{align} \label{eq: block_grad_bound_first_term_cond}
    &\ \mathbb{E}\left[\left.\left\Vert\tilde{\mathbf{g}}_k^{x,+}-\tilde{\mathbf{g}}_k^{x}\right\Vert^2\right|\mathcal{F}_k^b\right]\nonumber\\
    \leq&\ 3\cdot \mathbb{E}\left[\left.\left\Vert\tilde{\mathbf{g}}_k^{x,+}-\sum_{p\in\mathcal{J}_k^x}\nabla_{\mathbf{x},p}f(\mathbf{z}_k^+)\right\Vert^2+\left\Vert\tilde{\mathbf{g}}_k^{x}-\sum_{i\in\mathcal{I}_k^x}\nabla_{\mathbf{x},i}f(\mathbf{z}_k)\right\Vert^2\right|\mathcal{F}_k^b\right]\nonumber\\
    &\ +3\cdot \mathbb{E}\left[\left.\left\Vert\sum_{i\in\mathcal{I}_k^x}\nabla_{\mathbf{x},i}f(\mathbf{z}_k)-\sum_{p\in\mathcal{J}_k^x}\nabla_{\mathbf{x},p}f(\mathbf{z}_k^+)\right\Vert^2\right| \mathcal{F}_k^b\right]\nonumber\\
    \leq&\ \frac{3L^2r_k^2\tau_x}{2}+6\cdot\mathbb{E}\left[\left.\left\Vert\sum_{i\in\mathcal{I}_k^x}\nabla_{\mathbf{x},i}f(\mathbf{z}_k)\right\Vert^2+\left\Vert\sum_{p\in\mathcal{J}_k^x}\nabla_{\mathbf{x},p}f(\mathbf{z}_k^+)\right\Vert^2\right| \mathcal{F}_k^b\right]\nonumber\\
    \leq&\ \frac{3L^2r_k^2\tau_x}{2}+\frac{6\tau_x}{d_x}\mathbb{E}\left[\left.\left\Vert\nabla_\mathbf{x} f(\mathbf{z}_k)\right\Vert^2+\left\Vert\nabla_\mathbf{x} f(\mathbf{z}_k^+)\right\Vert^2\right|\mathcal{F}_k^b\right]\nonumber\\
     \leq&\ \ \frac{3L^2r_k^2\tau_x}{2}+\frac{12\tau_xG^2}{d_x},
\end{align}
where the second inequality follows from \eqref{eq: coor_grad_est_bias}.
Substituting \eqref{eq: block_grad_bound_first_term} and \eqref{eq: block_grad_bound_first_term_cond} into \eqref{eq: block_grad_bound_ZOCEG_term} gives:
\begin{equation}
\begin{aligned}\label{eq: block_grad_bound_x_term}
&\ \mathbb{E}\left[\left.\left\langle\nabla _{\mathbf{x}}f(\mathbf{z}_k^+),\mathbf{x}_k^+-\mathbf{x}\right\rangle\right| \mathcal{F}_k^b\right]\\  
\leq&\ \frac{d_x}{\eta\tau_x}\mathbb{E}\left[\left.H_k^x(\mathbf{x})-H_{k+1}^x(\mathbf{x})   \right| \mathcal{F}_k^b\right]+\frac{Lr_k\sqrt{d_x}\tilde{D}}{2} +\frac{3\eta L^2r_k^2d_x}{2}+12\eta G^2.
\end{aligned}  
\end{equation}
Similarly, for the dual variable, there holds:
\begin{equation}
\begin{aligned}\label{eq: block_grad_bound_y_term}
&\ \mathbb{E}\left[\left.-\left\langle\nabla _{\mathbf{y}}f(\mathbf{z}_k^+),\mathbf{y}_k^+-\mathbf{y}\right\rangle\right| \mathcal{F}_k^b\right]\\  
\leq&\ \frac{d_y}{\eta\tau_y}\mathbb{E}\left[\left. H_k^y(\mathbf{y})-H_{k+1}^y(\mathbf{y})  \right| \mathcal{F}_k^b\right]+\frac{Lr_k\sqrt{d_y}\tilde{D}}{2}   +\frac{3\eta L^2r_k^2d_y}{2}+12\eta G^2.
\end{aligned}  
\end{equation}
Combining \eqref{eq: block_grad_bound_x_term} and \eqref{eq: block_grad_bound_y_term} gives:
\begin{align*}   &\ \mathbb{E}\left[\left.\left\langle F(\mathbf{z}_k^+),\mathbf{z}_k^+-\mathbf{z}\right\rangle\right| \mathcal{F}_k^b\right]\\  
\leq&\ \frac{d_x}{\eta\tau_x}\mathbb{E}\left[\left. H_k^x(\mathbf{z})-H_{k+1}^x(\mathbf{z})   \right| \mathcal{F}_k^b\right]+\frac{d_y}{\eta\tau_y}\mathbb{E}\left[\left.H_k^y(\mathbf{z})-H_{k+1}^y(\mathbf{z})   \right| \mathcal{F}_k^b\right]+\frac{Lr_k(\sqrt{d_x}+\sqrt{d_y})\tilde{D}}{2}\nonumber\\
&\ +\frac{3(d_x+d_y)}{2}\eta L^2r_k^2+24\eta G^2.
\end{align*}
Taking the total expectation of and calculate the telescoping sum for $k=0,1,\cdots,K-1$, we can bound $\sum_{k=0}^{K-1}\mathbb{E}\left[\left\langle F(\mathbf{z}_k^+),\mathbf{z}_k^+-\mathbf{z}\right\rangle\right]$ as follows:
\begin{align}
\label{eq: bound_ZOBCEG_telescope}
&\ \sum_{k=0}^{K-1}\mathbb{E}\left[\left\langle F(\mathbf{z}_k^+),\mathbf{z}_k^+-\mathbf{z}\right\rangle\right]\nonumber\\
\leq &\ \frac{d_x}{\eta\tau_x}\mathbb{E}\left[H_0^x(\mathbf{z})-H_{K}^x(\mathbf{z})  \right]+\frac{d_y}{\eta\tau_y}\mathbb{E}\left[H_0^y(\mathbf{z})-H_{K}^y(\mathbf{z})  \right]+\frac{LM_1(\sqrt{d_x}+\sqrt{d_y})\tilde{D}}{2}\nonumber\\
&\ +\frac{3(d_x+d_y)}{2}\eta L^2 M_2+24\eta KG^2\nonumber\\
\leq &\ \frac{d_x\tilde{D}^2}{2\eta\tau_x}+\frac{d_y\tilde{D}^2}{2\eta\tau_y}+\frac{LM_1(\sqrt{d_x}+\sqrt{d_y})\tilde{D}}{2}+\frac{3(d_x+d_y)}{2}\eta L^2 M_2+24\eta KG^2, 
\end{align}
where the last inequality follows from the boundedness of $\mathcal{Z}$.
Substituting \eqref{eq: bound_ZOBCEG_telescope} into Lemma~\ref{lem: duality gap} finishes the convergence analysis of ZOBCEG.

\end{document}